\documentclass[12pt]{article}

\newtheorem{theorem}{Theorem}[section]
\newtheorem{lemma}[theorem]{Lemma}
\newtheorem{proposition}[theorem]{Proposition}
\newtheorem{corollary}[theorem]{Corollary}
\newtheorem{claim}[theorem]{Claim}

\newtheorem{example}[theorem]{Example}

\newtheorem{remark}[theorem]{Remark}

\newenvironment{poof}{\textit{Proof:  }}{
~\hfill\rule{2mm}{3mm}\vspace{.2in}}

\usepackage{amssymb,amsmath,tabularx,graphicx,float,placeins}
\usepackage{tikz}
\usetikzlibrary[calc,intersections,through,backgrounds,arrows,decorations.pathmorphing]

\def\Sfrak{\mathfrak{S}}

\def\Acal{\mathcal{A}}
\def\Bcal{\mathcal{B}}

\def\Fcal{\mathcal{F}}
\def\Jcal{\mathcal{J}}
\def\Ocal{\mathcal{O}}
\def\Pcal{\mathcal{P}}

\def\Cbb{\mathbb{C}}
\def\Nbb{\mathbb{N}}
\def\Rbb{\mathbb{R}}
\def\Zbb{\mathbb{Z}}

\def\ds{\displaystyle}

\def\ov{\overline}
\def\pa{\partial}
\def\pr{\prime}
\def\prpr{\prime\prime}
\def\ra{\rightarrow}
\def\Ra{\Rightarrow}
\def\setm{\setminus}
\def\ss{\scriptstyle}

\def\wtil{\widetilde}

\DeclareMathOperator{\Bic}{Bic}  
\DeclareMathOperator{\Camb}{Camb}  
\DeclareMathOperator{\Cell}{Cell}  
\DeclareMathOperator{\Con}{Con}  
\DeclareMathOperator{\con}{con}  
\DeclareMathOperator{\Cov}{Cov}  
\DeclareMathOperator{\GT}{GT}  
\DeclareMathOperator{\init}{init}  
\DeclareMathOperator{\lk}{lk}  
\DeclareMathOperator{\st}{st}  
\DeclareMathOperator{\term}{term}  
\DeclareMathOperator{\tr}{tr}

\begin{document}

\title{Lattice structure of Grid-Tamari orders}
\author{Thomas McConville}

\maketitle

\begin{abstract}
The Tamari order is a central object in algebraic combinatorics and many other areas.  Defined as the transitive closure of an associativity law, the Tamari order possesses a surprisingly rich structure:  it is a congruence-uniform lattice.  We consider a larger class of posets, the Grid-Tamari orders, which arise as an ordering on the facets of the non-kissing complex introduced by Pylyavskyy, Petersen, and Speyer.  In addition to Tamari orders, some interesting examples of Grid-Tamari orders include the Type A Cambrian lattices and Grassmann-Tamari orders.  We prove that the Grid-Tamari orders are congruence-uniform lattices, which resolves a conjecture of Santos, Stump, and Welker.  Towards this goal, we define a closure operator on sets of paths in a square grid, and prove that the biclosed sets of paths, ordered by inclusion, form a congruence-uniform lattice.  We then prove that the Grid-Tamari order is a quotient lattice of the corresponding lattice of biclosed sets.
\end{abstract}

\section{Introduction}

The Tamari lattice is a poset of proper bracketings of a word, with covering relations defined by the associativity law.  Tamari lattices and their generalizations have appeared in many parts of the literature.  We recommend the book \cite{tamari:festschrift} for an introduction to many recent developments on these posets.

We consider a new generalization of the Tamari lattice, the \emph{Grassmann-Tamari order}, introduced by Santos, Stump, and Welker \cite{santos.stump.welker:noncrossing}.  One of the conjectures they pose is that these posets are lattices.  We give an affirmative answer to this conjecture, and show that some of the very good lattice properties of Tamari lattices hold in this larger family of posets; see Theorem \ref{thm_main} for a precise statement.

The Grassmann-Tamari order $\GT_{k,n}$ is a partial order on the maximal ``non-crossing'' subsets of $\binom{[n]}{k}$, the $k$-element subsets of $\{1,\ldots,n\}$.  Two sets $I,J\in\binom{[n]}{k}$ are \emph{crossing} if $i_t<j_t<i_{t+1}<j_{t+1}$ for some $t$ where $I-J=\{i_1<\cdots<i_l\}$ and $J-I=\{j_1<\cdots<j_l\}$.  The sets $I,J$ are \emph{non-crossing} otherwise.  For example, $\{1,4,5\}$ and $\{2,3,6\}$ are non-crossing, whereas $\{1,4,5\}$ and $\{2,4,6\}$ are crossing.  The \emph{non-crossing complex} $\Delta^{NC}_{k,n}$ is the collection of all pairwise non-crossing subsets of $\binom{[n]}{k}$.

For $l\geq 1$, let $C_l$ be a chain poset with $l$ elements.  The complex $\Delta^{NC}_{k,n}$ may be realized as a regular, unimodular, Gorenstein triangulation of the order polytope $\Ocal_{k,n}$ on $C_k\times C_{n-k}$; i.e., the polytope in $\Rbb^{k(n-k)}$ defined by the inequalities $0\leq x_{i,j}\leq 1,\ x_{i,j}\leq x_{i+1,j}$, and $x_{i,j}\leq x_{i,j+1}$ for $1\leq i\leq k,\ 1\leq j\leq n-k$ (\cite[Theorem 8.1]{petersen.pylyavskyy.speyer:noncrossing} or \cite[Theorem 1.7]{santos.stump.welker:noncrossing}).  This triangulation of $\Ocal_{k,n}$ is distinct from the equatorial triangulation defined in \cite{reiner.welker:charney}, which is not flag in general.  As a consequence of this geometric realization, after removing cone points, $\Delta^{NC}_{k,n}$ is a pure, thin complex of dimension $(k-1)(n-k-1)-1$.  Moreover, there exists a simple polytope, the \emph{Grassmann-associahedron}, with facial structure anti-isomorphic to $\Delta^{NC}_{k,n}$.  As a flag, simplicial polytope, one may expect that the dual Grassmann-associahedron may be constructed by a sequence of suspensions and edge-stellations, which we prove in Section \ref{sec_polytope}.

Any triangulation of $\Ocal_{k,n}$ naturally gives rise to a monomial basis for the coordinate ring of the Grassmannian, the $\Cbb$-algebra generated by the maximal minors of a $k\times n$ matrix of indeterminates $(x_{ij})$ \cite{petersen.pylyavskyy.speyer:noncrossing}.  Namely, a monomial $\prod_1^r x_{I_j}$ is in the basis if $\{I_1,\ldots,I_r\}$ is a face of the triangulation.  The classical standard basis for this algebra is indexed by semistandard Young tableaux.  The columns of a semistandard Young tableaux satisfy a compatibility condition that resembles a non-nesting analogue of the non-crossing condition defined above.  Thus these two bases may be viewed as ``opposite'' in some sense; see \cite[Remark 4.7]{santos.stump.welker:noncrossing}.  One may hope to develop a straightening law for these monomials, though we do not pursue this here.

Let $\Jcal$ be the set of order ideals of $C_k\times C_{n-k}$.  The \emph{Hibi ideal} is the ideal generated by $\{x_Ix_J-x_{I\cap J}x_{I\cup J}:\ I,J\in\Jcal)$ in the polynomial ring on $\{x_I:\ I\in\Jcal\}$.  By results of \cite{sturmfels:grobner}, regular unimodular triangulations of $\Ocal_{k,n}$ are in bijection with squarefree monomial initial ideals of the Hibi ideal.  As observed in the introduction of \cite{santos.stump.welker:noncrossing}, the triangulation induced by $\Delta^{NC}_{k,n}$ corresponds to a particularly nice initial ideal.  We refer to the survey \cite[Section 6]{conca.hosten.thomas:nice} for more background on Hibi ideals.

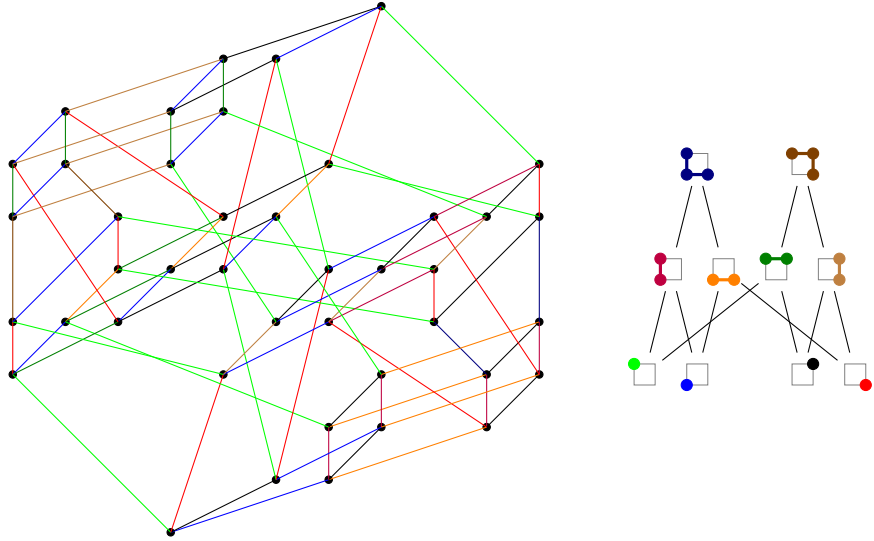
\begin{figure}
\begin{centering}

\begin{tikzpicture}[scale=.7]
\colorlet{e1}{blue}
\colorlet{e2}{black}
\colorlet{e3}{green}
\colorlet{e4}{red}
\colorlet{e5}{purple}
\colorlet{e6}{orange}
\colorlet{e7}{green!50!black}
\colorlet{e8}{brown}
\colorlet{e9}{blue!50!black}
\colorlet{e10}{orange!50!black}

\begin{scope}

\coordinate (v1) at (-1,0);
\coordinate (v2) at (2,1);
\coordinate (v3) at (2,2); 
\coordinate (v4) at (5,2); 
\coordinate (v5) at (1,1); 
\coordinate (v6) at (-4,3); 
\coordinate (v7) at (0,3); 
\coordinate (v8) at (5,3); 
\coordinate (v9) at (3,2);
\coordinate (v10) at (-4,4);
\coordinate (v11) at (-3,4);
\coordinate (v12) at (-2,4);
\coordinate (v13) at (4,4);
\coordinate (v14) at (2,4); 
\coordinate (v15) at (1,4); 
\coordinate (v16) at (3,3); 
\coordinate (v17) at (6,3); 
\coordinate (v18) at (-2,5);
\coordinate (v19) at (-1,5);
\coordinate (v20) at (0,5);
\coordinate (v21) at (4,5); 
\coordinate (v22) at (3,5); 
\coordinate (v23) at (6,4); 
\coordinate (v24) at (2,5); 
\coordinate (v25) at (-4,6);
\coordinate (v26) at (-2,6);
\coordinate (v27) at (0,6);
\coordinate (v28) at (1,6); 
\coordinate (v29) at (5,6); 
\coordinate (v30) at (6,6); 
\coordinate (v31) at (4,6); 
\coordinate (v32) at (-4,7);
\coordinate (v33) at (-3,7);
\coordinate (v34) at (-1,7);
\coordinate (v35) at (2,7); 
\coordinate (v36) at (6,7); 
\coordinate (v37) at (-3,8);
\coordinate (v38) at (-1,8);
\coordinate (v39) at (0,8); 
\coordinate (v40) at (0,9);
\coordinate (v41) at (1,9); 
\coordinate (v42) at (3,10);

\foreach \v in {(v1),(v2),(v3),(v4),(v5),(v6),(v7),(v8),(v9),(v10),(v11),(v12),(v13),(v14),(v15),(v16),(v17),(v18),(v19),(v20),(v21),(v22),(v23),(v24),(v25),(v26),(v27),(v28),(v29),(v30),(v31),(v32),(v33),(v34),(v35),(v36),(v37),(v38),(v39),(v40),(v41),(v42)}
  {\filldraw \v circle(2pt); }

\draw[color=e1] (v1) -- (v2);
\draw[color=e2] (v1) -- (v5);
\draw[color=e3] (v1) -- (v6);
\draw[color=e4] (v1) -- (v7);
\draw[color=e5] (v2) -- (v3);
\draw[color=e6] (v2) -- (v4);
\draw[color=e2] (v2) -- (v9);
\draw[color=e6] (v3) -- (v8);
\draw[color=e3] (v3) -- (v11);
\draw[color=e2] (v3) -- (v16);
\draw[color=e5] (v4) -- (v8);
\draw[color=e4] (v4) -- (v14);
\draw[color=e2] (v4) -- (v17);
\draw[color=e1] (v5) -- (v9);
\draw[color=e3] (v5) -- (v20);
\draw[color=e4] (v5) -- (v24);
\draw[color=e4] (v6) -- (v10);
\draw[color=e1] (v6) -- (v11);
\draw[color=e7] (v6) -- (v12);
\draw[color=e3] (v7) -- (v10);
\draw[color=e1] (v7) -- (v14);
\draw[color=e8] (v7) -- (v15);
\draw[color=e9] (v8) -- (v13);
\draw[color=e2] (v8) -- (v23);
\draw[color=e5] (v9) -- (v16);
\draw[color=e6] (v9) -- (v17);
\draw[color=e10] (v10) -- (v25);
\draw[color=e1] (v10) -- (v26);
\draw[color=e6] (v11) -- (v18);
\draw[color=e7] (v11) -- (v19);
\draw[color=e1] (v12) -- (v19);
\draw[color=e2] (v12) -- (v20);
\draw[color=e4] (v12) -- (v32);
\draw[color=e3] (v13) -- (v18);
\draw[color=e4] (v13) -- (v21);
\draw[color=e2] (v13) -- (v30);
\draw[color=e5] (v14) -- (v21);
\draw[color=e8] (v14) -- (v22);
\draw[color=e1] (v15) -- (v22);
\draw[color=e2] (v15) -- (v24);
\draw[color=e3] (v15) -- (v34);
\draw[color=e6] (v16) -- (v23);
\draw[color=e3] (v16) -- (v28);
\draw[color=e5] (v17) -- (v23);
\draw[color=e4] (v17) -- (v31);
\draw[color=e4] (v18) -- (v26);
\draw[color=e7] (v18) -- (v27);
\draw[color=e6] (v19) -- (v27);
\draw[color=e2] (v19) -- (v28);
\draw[color=e1] (v20) -- (v28);
\draw[color=e4] (v20) -- (v41);
\draw[color=e3] (v21) -- (v26);
\draw[color=e8] (v21) -- (v29);
\draw[color=e5] (v22) -- (v29);
\draw[color=e2] (v22) -- (v31);
\draw[color=e9] (v23) -- (v30);
\draw[color=e1] (v24) -- (v31);
\draw[color=e3] (v24) -- (v41);
\draw[color=e7] (v25) -- (v32);
\draw[color=e1] (v25) -- (v33);
\draw[color=e8] (v25) -- (v34);
\draw[color=e10] (v26) -- (v33);
\draw[color=e2] (v27) -- (v35);
\draw[color=e4] (v27) -- (v37);
\draw[color=e6] (v28) -- (v35);
\draw[color=e2] (v29) -- (v36);
\draw[color=e3] (v29) -- (v39);
\draw[color=e3] (v30) -- (v35);
\draw[color=e4] (v30) -- (v36);
\draw[color=e5] (v31) -- (v36);
\draw[color=e1] (v32) -- (v37);
\draw[color=e8] (v32) -- (v38);
\draw[color=e7] (v33) -- (v37);
\draw[color=e8] (v33) -- (v39);
\draw[color=e7] (v34) -- (v38);
\draw[color=e1] (v34) -- (v39);
\draw[color=e4] (v35) -- (v42);
\draw[color=e3] (v36) -- (v42);
\draw[color=e8] (v37) -- (v40);
\draw[color=e1] (v38) -- (v40);
\draw[color=e2] (v38) -- (v41);
\draw[color=e7] (v39) -- (v40);
\draw[color=e2] (v40) -- (v42);
\draw[color=e1] (v41) -- (v42);

\end{scope}

\begin{scope}[xshift=10cm,yshift=3cm]
\scriptsize;
\draw[gray, very thin] (-2.2,-.2) rectangle (-1.8,.2);
\filldraw[color=e3] (-2.2,.2) circle(3pt);

\draw[gray, very thin] (-1.2,-.2) rectangle (-.8,.2);
\filldraw[color=e1] (-1.2,-.2) circle(3pt);

\draw[gray, very thin] (.8,-.2) rectangle (1.2,.2);
\filldraw[color=e2] (1.2,.2) circle(3pt);

\draw[gray, very thin] (1.8,-.2) rectangle (2.2,.2);
\filldraw[color=e4] (2.2,-.2) circle(3pt);

\draw[gray, very thin] (-1.7,1.8) rectangle (-1.3,2.2);
\filldraw[color=e5] (-1.7,2.2) circle(3pt)
                    (-1.7,1.8) circle(3pt);
\draw[color=e5, very thick] (-1.7,2.2) -- (-1.7,1.8);

\draw[gray, very thin] (-.7,1.8) rectangle (-.3,2.2);
\filldraw[color=e6] (-.7,1.8) circle(3pt)
                    (-.3,1.8) circle(3pt);
\draw[color=e6, very thick] (-.7,1.8) -- (-.3,1.8);

\draw[gray, very thin] (.3,1.8) rectangle (.7,2.2);
\filldraw[color=e7] (.3,2.2) circle(3pt)
                    (.7,2.2) circle(3pt);
\draw[color=e7, very thick] (.3,2.2) -- (.7,2.2);

\draw[gray, very thin] (1.3,1.8) rectangle (1.7,2.2);
\filldraw[color=e8] (1.7,1.8) circle(3pt)
                    (1.7,2.2) circle(3pt);
\draw[color=e8, very thick] (1.7,1.8) -- (1.7,2.2);

\draw[gray, very thin] (-1.2,3.8) rectangle (-.8,4.2);
\filldraw[color=e9] (-1.2,4.2) circle(3pt)
                     (-1.2,3.8) circle(3pt)
                     (-.8,3.8) circle(3pt);
\draw[color=e9, very thick] (-1.2,4.2) -- (-1.2,3.8) -- (-.8,3.8);

\draw[gray, very thin] (.8,3.8) rectangle (1.2,4.2);
\filldraw[color=e10] (.8,4.2) circle(3pt)
                     (1.2,4.2) circle(3pt)
                     (1.2,3.8) circle(3pt);
\draw[color=e10, very thick] (.8,4.2) -- (1.2,4.2) -- (1.2,3.8);

\end{scope}

\begin{scope}[xshift=10cm,yshift=3cm,xscale=.5,yscale=2]

\draw[shorten <=.3cm, shorten >=.3cm] (-4,0) -- (-3,1);
\draw[shorten <=.3cm, shorten >=.3cm] (-3,1) -- (-2,0);
\draw[shorten <=.3cm, shorten >=.3cm] (-2,0) -- (-1,1);
\draw[shorten <=.3cm, shorten >=.3cm] (-1,1) -- (4,0);
\draw[shorten <=.3cm, shorten >=.3cm] (4,0) -- (3,1);
\draw[shorten <=.3cm, shorten >=.3cm] (3,1) -- (2,0);
\draw[shorten <=.3cm, shorten >=.3cm] (2,0) -- (1,1);
\draw[shorten <=.3cm, shorten >=.3cm] (1,1) -- (-4,0);
\draw[shorten <=.3cm, shorten >=.3cm] (-3,1) -- (-2,2);
\draw[shorten <=.3cm, shorten >=.3cm] (-2,2) -- (-1,1);
\draw[shorten <=.3cm, shorten >=.3cm] (3,1) -- (2,2);
\draw[shorten <=.3cm, shorten >=.3cm] (2,2) -- (1,1);

\end{scope}

\end{tikzpicture}
\caption{\label{fig_P36}(left) $\GT_{3,6}$ (right) $J(\Con(\GT_{3,6}))^*$}
\end{centering}
\end{figure}

There is a natural orientation on the dual graph of $\Delta^{NC}_{k,n}$.  If two facets $F_1=F\cup\{I\},F_2=F\cup\{J\}$ are adjacent, then there is a unique index $t$ for which $i_t<j_t<i_{t+1}<j_{t+1}$ where $I-J=\{i_1<\cdots<i_l\}$ and $J-I=\{j_1<\cdots<j_l\}$.  We orient the edge $F_1\ra F_2$ if the pair $\{i_t,i_{t+1}\}$ is lexicographically smaller than $\{j_t,j_{t+1}\}$.  For example, $\{145,146,236,245\}$ and $\{146,236,245,246\}$ are adjacent facets of $\Delta^{NC}_{3,6}$ with orientation $\{145,146,236,245\}\ra\{146,236,245,246\}$ since $145$ and $246$ cross at $15$ and $26$.  Defined by Santos, Stump, and Welker in \cite{santos.stump.welker:noncrossing}, the \emph{Grassmann-Tamari order} $\GT_{k,n}$ is the transitive closure of this relation.  The smallest Grassmann-Tamari order not isomorphic to a Cambrian lattice is drawn in Figure \ref{fig_P36}.

\begin{figure}
\begin{centering}
\begin{tikzpicture}

\begin{scope}
\draw[step=1cm] (0,0) grid (3,3);
\draw[decorate,decoration={coil,segment length=4pt},color=red,xshift=-2mm,yshift=-2mm] (.2,2) -- (2,2) -- (2,.2);
\draw[decorate,decoration={coil,segment length=4pt},color=green,xshift=2mm,yshift=2mm] (1,2.8) -- (1,2) -- (2,2) -- (2,1) -- (2.8,1);
\draw[ultra thick] (2,1) -- (2,2) -- (1,2);
\filldraw (2,1) circle(1mm);
\filldraw (1,2) circle(1mm);
\draw (.8,2.2) node{$\ss v$};
\draw (2.2,.8) node{$\ss v^{\pr}$};
\end{scope}

\begin{scope}[xshift=6cm]

\draw[step=1cm] (0,0) grid (3,3);
\draw[decorate,decoration={coil,segment length=4pt},color=red] (0,2) -- (2,2) -- (2,0);
\draw[decorate,decoration=zigzag,color=blue,xshift=2mm,yshift=2mm] (-.2,2) -- (2,2) -- (2,1) -- (2.8,1);
\draw[decorate,decoration={coil,segment length=4pt},color=green] (1,3) -- (1,1) -- (3,1);
\draw[decorate,decoration=saw,color=orange,xshift=-2mm,yshift=-2mm] (1,3.2) -- (1,2) -- (2,2) -- (2,0.2);

\end{scope}

\end{tikzpicture}
\caption{\scriptsize\label{fig_kissing}(left) Two paths kissing along the indicated segment from $v$ to $v^{\pr}$.  The paths correspond to the sets $145$ and $246$, which are crossing. (right) A maximal family of non-kissing paths excluding horizontal and vertical paths.}
\end{centering}
\end{figure}
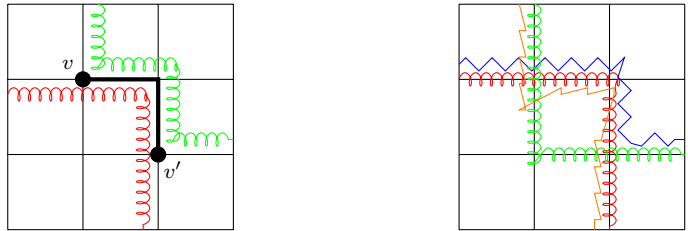

The non-crossing condition translates to a \emph{non-kissing} condition on paths via the standard bijection between $k$-subsets of $[n]$ and paths in a $k\times(n-k)$ rectangle with South and East steps.  For example, the set $\{1,4,5\}$ corresponds to the path from the NW-corner to the SE-corner of the rectangle such that the first, fourth, and fifth steps are to the South, while the others are to the East.  Two paths $p_1,p_2$ in the plane are \emph{kissing} if they agree on some subpath between vertices $v$ and $v^{\pr}$ such that
\begin{enumerate}
\item $p_1$ enters $v$ from the West and leaves $v^{\pr}$ to the South, and 
\item $p_2$ enters $v$ from the North and leaves $v^{\pr}$ to the East.
\end{enumerate}
An example of two kissing paths is given in Figure \ref{fig_kissing}.  The non-kissing complex $\Delta^{NK}(\lambda)$ associated to a (possibly not rectangular) shape $\lambda$ is the collection of pairwise non-kissing paths supported by $\lambda$.  A poset $\GT(\lambda)$ analogous to the Grassmann-Tamari orders may be defined on the facets of this complex.  We call $\GT(\lambda)$ the \emph{Grid-Tamari order}; see Section \ref{sec_non-kissing}.

Our main result is

\begin{theorem}\label{thm_main}
For any shape $\lambda$, $\GT(\lambda)$ is a congruence-uniform lattice.
\end{theorem}

We recall congruence-uniformity and related lattice properties in Section \ref{sec_biclosed}.

To prove Theorem \ref{thm_main}, we express $\GT(\lambda)$ as a lattice quotient of a much simpler lattice.  Namely, we define a finite topological space whose clopen sets, which we call \emph{biclosed sets}, form a congruence-uniform lattice under inclusion.  Then we define a map from the collection of biclosed sets to facets of the non-kissing complex that carries this lattice structure.

The paper is organized as follows.  Some notation and basic results on lattices are introduced in Section \ref{sec_lattices}.  In Section \ref{sec_non-kissing}, we establish the purity and thinness of the non-kissing complex combinatorially, similar to the methodology employed in \cite[Section 2.2]{santos.stump.welker:noncrossing} for proving purity and thinness of the non-crossing complex.  We close the section by defining the orientation on the dual graph of the non-kissing complex whose transitive closure is a Grid-Tamari order.  We emphasize that this directed graph is acyclic as a \emph{consequence} of Theorem \ref{thm_main}.  A geometric proof of acyclicity in the non-crossing case appears in \cite{santos.stump.welker:noncrossing}.

The reduced non-kissing complex is the boundary complex of a simplicial polytope.  In Section \ref{sec_polytope}, we describe a way to construct these polytopes by a sequence of edge-stellations and suspensions.  It follows that the $\gamma$-vector of the non-kissing complex is the $f$-vector of a flag simplicial complex by results of \cite{aisbett:gamma} and \cite{buchstaber.volodin:combinatorial}.  Moreover, the Hasse diagram for the Grid-Tamari order is the 1-skeleton of the polar dual polytope.  These dual polytopes may be constructed by dual operations, namely ridge truncations and doublings.  We remark that although ridge truncations sometimes correspond to interval doublings, these two constructions do not match up in general.

We prove some general results on biclosed sets in Section \ref{sec_biclosed}.  Biclosed sets may be defined for any closure operator on a set, though the resulting poset of biclosed sets may not be interesting.  We provide some conditions on the closure that makes the poset of biclosed sets a congruence-uniform lattice in Theorem \ref{thm_closure_CN}.  In particular, these conditions are satisfied by the convex closure on the positive roots of a finite root system.

In Section \ref{sec_segments}, we introduce a poset of biclosed subsets of segments in a shape $\lambda$.  A collection of segments between two interior vertices of $\lambda$ is \emph{closed} if a segment $s$ is in $\lambda$ whenever there exists a partition of $s$ into two subpaths that both lie in $\lambda$.  We show that this closure satisfies the hypotheses given in Section \ref{sec_biclosed}, so its poset of biclosed sets is a congruence-uniform lattice.

A special lattice congruence on the lattice of biclosed sets of segments is presented in Section \ref{sec_quotient}.  In Section \ref{sec_tamari}, we define a map $\eta$ from biclosed sets of segments to the facets of the non-kissing complex, and show that the fibers of $\eta$ are precisely the equivalence classes of this lattice congruence.  We then deduce Theorem \ref{thm_main} by comparing the order induced by $\eta$ with the Grid-Tamari order.

Some other interesting lattice quotients of the weak order called \emph{Cambrian lattices} were introduced by Reading in \cite{reading:cambrian_lattices}.  In Section \ref{sec_cambrian}, we prove that the type $A$ Cambrian lattices are examples of Grid-Tamari orders for double ribbon shapes.  We prove this isomorphism using Reading's description of type $A$ Cambrian lattices as a poset of triangulations of a polygon.

\section{Lattices}\label{sec_lattices}

In this section, we set up some notation for lattices, mostly following \cite{gratzer:congruences}.

A \emph{lattice} is a partially ordered set (\emph{poset}) for which any two elements $x,y$ have a least upper bound $x\vee y$ and a greatest lower bound $x\wedge y$.  If $x<y$, we say $y$ \emph{covers} $x$ if there does not exist $z$ such that $x<z<y$.  If $P$ is a poset with an element $x$ such that $x\leq y$ for all $y\in P$, then $x$ is the \emph{bottom} element of $P$, typically denoted $\hat{0}$.  Dually, the \emph{top} element of $P$ is denoted $\hat{1}$.  An \emph{order ideal} $X$ of a poset $P$ is a subset of $P$ such that if $x\leq y$ and $y\in X$ then $x\in X$.  We let $\Ocal(P)$ denote the set of order ideals of $P$.  The \emph{dual poset} $P^*$ has the same underlying set as $P$ where $x\leq_{P^*}y$ if and only if $y\leq_P x$.

The following lemma is frequently used to prove that a poset is a lattice, see e.g. \cite{bjorner.edelman.ziegler:lattice}, \cite{hersh.meszaros:sb}, \cite{jambu.paris:combinatorics}, \cite{mcconville:crosscut}.

\begin{lemma}[\cite{bjorner.edelman.ziegler:lattice} Lemma 2.1]\label{lem_local_lattice}
Let $P$ be a finite poset with $\hat{0}$ and $\hat{1}$.  If $x\vee y$ exists for $x,y,z\in P$ such that $x$ and $y$ both cover $z$, then $P$ is a lattice.
\end{lemma}

In Lemma \ref{lem_local_map_join}, we describe a similar result for maps between lattices.

\begin{lemma}\label{lem_local_map_join}
Let $f:L\ra L^{\pr}$ be an order-preserving map between finite lattices $L$ and $L^{\pr}$.
\begin{enumerate}
\item\label{lem_lmj_1} Suppose $f(x\vee y)=f(x)\vee f(y)$ for $x,y,z\in L$ such that $x$ and $y$ both cover $z$.  Then $f(a\vee b)=f(a)\vee f(b)$ for all $a,b\in L$.
\item\label{lem_lmj_2} Suppose $f(x)=f(y)$ implies $f(x\vee y)=f(x)$ for $x,y,z\in L$ such that $x$ and $y$ both cover $z$.  If $f$ preserves meets, then $f(a)=f(b)$ implies $f(a\vee b)=f(a)$ for all $a,b\in L$.
\end{enumerate}
\end{lemma}

\begin{poof}
For $x\in L$, define the depth of $x$ to be the length of the longest chain from $x$ to $\hat{1}$.  We prove both statements by induction on depth.

(\ref{lem_lmj_1}): Let $a,b\in L$.  Assume $f(a^{\pr}\vee b^{\pr})=f(a^{\pr})\vee f(b^{\pr})$ whenever $a\wedge b<a^{\pr}\wedge b^{\pr}$.  If $a\leq b$ then $f(a\vee b)=f(b)=f(a)\vee f(b)$ since $f$ is order-preserving.

Assume $a$ and $b$ are incomparable, and let $x$ and $y$ cover $a\wedge b$ such that $x\leq a$ and $y\leq b$.  Then $x\neq y$ and $f(x\vee y)=f(x)\vee f(y)$ by assumption.  Since $x\leq a\wedge(x\vee y)$ holds, we have $f(a)\vee f(x\vee y)=f(a\vee x\vee y)=f(a\vee y)$ by induction.  Similarly, $f(b)\vee f(x\vee y)=f(b\vee x\vee y)=f(b\vee x)$ holds.  Since $x\vee y\leq(a\vee y)\wedge(b\vee x)$, we deduce
\begin{align*}
f(a\vee b)=f(a\vee y\vee b\vee x)=f(a\vee y)\vee f(b\vee x) &= f(a)\vee f(x\vee y)\vee f(b)\\
&= f(a)\vee f(b)\vee f(x)\vee f(y)\\
&= f(a)\vee f(b).
\end{align*}

(\ref{lem_lmj_2}): Assume $f$ preserves meets.  Let $a,b\in L$ such that $f(a)=f(b)$, and set $w=f(a)$.  If $a\leq b$, then $f(a\vee b)=f(b)=f(a)$ holds.

Assume $a$ and $b$ are incomparable, and let $x$ and $y$ cover $a\wedge b$ such that $x\leq a$ and $y\leq b$.  Since $f(a)=f(b)=w$ and $f$ preserves meets, we have $f(a\wedge b)=w$.  As $f$ is order-preserving, this implies $f(x)=w=f(y)$.  In particular, $f(x\vee y)=w$ by assumption.  As before, we deduce that $f(a\vee(x\vee y))=w$ and $f(b\vee(x\vee y))=w$ by the induction hypothesis.  Applying the induction hypothesis again, we deduce $f(a\vee b)=w$.
\end{poof}

An equivalence relation $\Theta$ on a lattice $L$ is a \emph{lattice congruence} if $x\equiv y\mod\Theta$ implies $x\vee z\equiv y\vee z\mod\Theta$ and $x\wedge z\equiv y\wedge z\mod\Theta$ for $x,y,z\in L$.  The set of equivalence classes $L/\Theta$ of a lattice congruence forms a lattice where $[x]\vee[y]=[x\vee y]$ and $[x]\wedge[y]=[x\wedge y]$ for $x,y\in L$.  We say $L/\Theta$ is a \emph{quotient lattice} of $L$, and the natural map $L\mapsto L/\Theta$ is a \emph{lattice quotient map}.  The following characterization of lattice congruences is well-known.

\begin{proposition}\label{prop_congruence}
Let $\Theta$ be an equivalence relation on a finite lattice $L$.  If
\begin{enumerate}
\item the equivalence classes of $\Theta$ are all closed intervals of $L$, and
\item the maps $\pi^{\uparrow}$ and $\pi_{\downarrow}$ taking an element of $L$ to the largest (respectively, smallest) element of its equivalence class are both order-preserving,
\end{enumerate}
then $\Theta$ is a lattice congruence.
\end{proposition}

An element $j$ of a lattice $L$ is \emph{join-irreducible} if for $x,y\in L$ such that $j=x\vee y$, either $j=x$ or $j=y$.  If $L$ is finite, $j$ is join-irreducible exactly when it covers a unique element, which we call $j_*$.  A \emph{meet-irreducible} element $m$ is defined dually and is covered by a unique element $m^*$.  We let $J(L)$ and $M(L)$ denote the sets of join-irreducible and meet-irreducible elements of $L$, respectively.

Given a lattice $L$, its set of lattice congruences $\Con(L)$ forms a distributive lattice under refinement order.  Hence when $L$ is finite, $\Con(L)$ is isomorphic to $\Ocal(J(\Con(L)))$.  If $y$ covers $x$, we write $\con(x,y)$ for the minimal lattice congruence in which $x\equiv y\ (\con(x,y))$ holds.

For any finite lattice $L$ with lattice congruence $\Theta$, we have
$$\Theta=\bigvee_{\substack{j\in J(L)\\ j\equiv j_*\mod\Theta}}\con(j_*,j).$$
Hence, the join-irreducible congruences are always of the form $\con(j_*,j)$ for some $j\in J(L)$.  A finite lattice $L$ is \emph{congruence-uniform} (or \emph{bounded}) if
\begin{itemize}
\item the map $j\mapsto\con(j_*,j)$ is a bijection from $J(L)$ to $J(\Con(L))$, and
\item the map $m\mapsto\con(m,m^*)$ is a bijection from $M(L)$ to $M(\Con(L))$.
\end{itemize}
Alternatively, finite congruence-uniform lattices may be characterized as homomorphic images of free lattices with bounded fibers or as lattices constructible from the one-element lattice by a sequence of interval doublings \cite{day:congruence}; see Section \ref{subsec_CN} for this construction.

\section{Non-kissing complexes}\label{sec_non-kissing}

Let $\lambda$ be a finite induced subgraph of the $\Zbb\times\Zbb$ square grid.  We refer to such a graph as a \emph{shape}.  A vertex $v$ is \emph{interior} if $\lambda$ contains the $2\times 2$ grid centered at $v$.  Any vertex of $\lambda$ that is not interior is called a \emph{boundary vertex}.  We say $v$ is a \emph{SE-corner} if the vertices one step South or East of $v$ are not in $\lambda$.  If $v$ is a vertex of $\lambda$, then $\lambda\setm v$ is the subgraph of $\lambda$ with $v$ removed.

A \emph{path} supported by $\lambda$ is a sequence of vertices $v_0,\ldots,v_t$ such that
\begin{itemize}
\item $v_0$ and $v_t$ are boundary vertices,
\item $v_1,\ldots,v_{t-1}$ are interior vertices, and
\item $v_i$ is one step South or East of $v_{i-1}$ for all $i$.
\end{itemize}

\begin{figure}
\begin{centering}
\begin{tikzpicture}

\draw[step=1cm, gray, very thin] (0,-2) grid (5,0);
\draw[step=1cm, gray, very thin] (1,-3) grid (4,-2);
\draw[step=1cm, gray, very thin] (2,-4) grid (4,-3);

\draw[decorate,decoration={coil,segment length=4pt},color=red] (1,-2) -- (3,-2) -- (3,-3) -- (4,-3);
\draw (1,-2) node[anchor=south east]{$w_0$};
\draw (2,-2) node[anchor=south east]{$w_1$};
\draw (3,-2) node[anchor=south west]{$w_2$};
\draw (3,-3) node[anchor=north east]{$w_3$};
\draw (4,-3) node[anchor=north west]{$w_4$};
\fill[black] (4,-4) circle(4pt);
\draw (4,-4) node[anchor=north west]{$v$};

\end{tikzpicture}
\caption{\scriptsize\label{fig_shape}A shape with a path $w_0,w_1,w_2,w_3,w_4$ and a SE-corner $v$.}
\end{centering}
\end{figure}
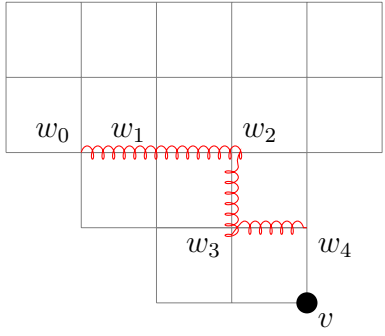

\begin{example}
For the path in Figure \ref{fig_shape}, the vertices $w_0$ and $w_4$ are boundary vertices, while $w_1,w_2,$ and $w_3$ are interior.  If $v$ is removed from $\lambda$, then $w_3$ becomes a boundary vertex.  The restriction of this path to $\lambda\setm v$ is the sequence $w_0,w_1,w_2,w_3$.
\end{example}

A path supported by $\lambda$ is called a \emph{segment} if its endpoints are also interior vertices.  If $s$ is a segment containing vertices $v$ and $v^{\pr}$, then $s[v,v^{\pr}]$ denotes the sub-segment of $s$ whose endpoints are $v$ and $v^{\pr}$.  The initial (terminal) vertex of a segment $s$ is denoted $s_{\init}$ ($s_{\term}$).  We abbreviate $s[s_{\init},v]$ and $s[v,s_{\term}]$ to $s[\cdot,v]$ and $s[v,\cdot]$, respectively.  A segment that only contains one vertex is called \emph{lazy}.  All other segments are \emph{non-lazy}.

Two paths $p_1,p_2$ are \emph{kissing} if they share vertices $v,v^{\pr}$ such that
\begin{itemize}
\item $p_1[v,v^{\pr}]=p_2[v,v^{\pr}]$,
\item $p_1$ enters $v$ from the West and leaves $v^{\pr}$ to the South, and
\item $p_2$ enters $v$ from the North and leaves $v^{\pr}$ to the East.
\end{itemize}

Otherwise $p_1$ and $p_2$ are \emph{non-kissing}.  The \emph{non-kissing complex} $\Delta^{NK}(\lambda)$ is the (flag) simplicial complex whose faces are collections of pairwise non-kissing paths supported by $\lambda$.  Let $\Fcal(\Delta^{NK}(\lambda))$ denote the set of \emph{facets}, the maximal faces of this complex.  As horizontal and vertical paths are non-kissing with any path, we define the \emph{reduced non-kissing complex} $\wtil{\Delta}^{NK}(\lambda)$ to be the deletion of all horizontal and vertical paths from $\Delta^{NK}(\lambda)$.

Although a pair of non-kissing paths may twist around each other several times, there is a natural way to totally order paths that contain a specific edge.  Let $e$ be an edge of $\lambda$.  If $p_1$ and $p_2$ are distinct non-kissing paths containing $e$, then they agree on some maximal segment $p_1[v,v^{\pr}]$ containing $e$.  Order $p_1\prec_e p_2$ if either $p_1$ enters $v$ from the North or $p_1$ leaves $v^{\pr}$ to the South.  A path $p\in F$ is the \emph{bottom path} (\emph{top path}) at an edge $e$ if $p$ is minimal (maximal) in $F$ with respect to $\prec_e$.

\begin{theorem}\label{thm_pure_thin}
Let $F$ be a facet of $\Delta^{NK}(\lambda)$.
\begin{enumerate}
\item\label{thm_pure_thin_1} The map $e\mapsto\max_{\prec_e}F$ is a bijection between vertical edges of $\lambda$ and non-horizontal paths in $F$.
\item\label{thm_pure_thin_2} Dually, the map $e\mapsto\min_{\prec_e}F$ is a bijection between horizontal edges of $\lambda$ and non-vertical paths in $F$.
\item\label{thm_pure_thin_3} For paths $p\in F$ with at least one turn, there exists a unique path $q$ distinct from $p$ such that $F-\{p\}\cup\{q\}$ is non-kissing.  Moreover, $p$ and $q$ kiss at a unique segment.
\end{enumerate}
\end{theorem}

\begin{poof}
For each of these statements, we proceed by induction on the size of $\lambda$.  Let $c$ be SE-corner of $\lambda$ and let $w$ be the point in $\Zbb\times\Zbb$ one step NW of $c$.  If $w$ is not an interior vertex of $\lambda$, then every path in $\lambda$ is supported by $S\setm c$, so the theorem holds by the inductive hypothesis.  Hence, we may assume that $w$ is an interior vertex of $\lambda$.

(\ref{thm_pure_thin_1}): We start by proving injectivity of the map.  Suppose there is a path $p\in F$ that is on top at two distinct vertical edges $e_1,e_2$.  Let $v$ be the southern vertex of $e_1$ and let $e$ be the edge west of $v$.  Let $p^{\pr}\in F$ be the bottom path at $e$.  Define a path $q$ supported by $\lambda$ where $q[\cdot,v]=p^{\pr}[\cdot,v]$ and $q[v,\cdot]=p[v,\cdot]$.  Since $p\prec_{e_2} q$, the path $q$ is not in $F$.  Let $t$ be the segment containing $v$ along which $p$ and $p^{\pr}$ agree.  Since $p$ and $p^{\pr}$ are non-kissing, $p$ leaves $t$ to the South and $p^{\pr}$ leaves to the East.

We claim that $q$ is non-kissing with every path in $F$, contradicting the maximality of $F$.  Indeed, if $q$ and $q^{\pr}$ are kissing for some $q^{\pr}\in F$, then they must kiss at a segment $s$ containing $v$ as $q^{\pr}$ is non-kissing with both $p$ and $p^{\pr}$.

If $v$ is the initial vertex of $s$, then $q^{\pr}$ must leave the terminal vertex of $s$ to the east while $q$ leaves to the south.  But this means $q^{\pr}$ enters $v$ from the North, which contradicts maximality of $p$ at $e_1$.

If $v$ is not the initial vertex of $s$, then $q^{\pr}$ contains $e$.  By the minimality of $p^{\pr}$ at $e$, $q^{\pr}$ must enter $s$ from the West and leave $s$ to the South.  If $t$ is a subsegment of $s$, then $q^{\pr}\prec_e p^{\pr}$, a contradiction.  If $t$ contains $s$, then $q^{\pr}$ and $p$ are kissing, a contradiction.

Next we verify surjectivity.  The restriction of paths in $F$ to $S\setm c$ defines a collection of non-kissing paths supported by $S\setm c$.  Let $e_1$ be the edge south of $w$ and $e_2$ the edge east of $w$.  It is straight-forward to check that if $p$ and $p^{\pr}$ are distinct paths on $\lambda$ with the same restriction to $S\setm c$, then $p$ must the top path at $e_1$ and $p^{\pr}$ the bottom path at $e_2$ (or vice versa).  Hence, the map applied to $F\setm c$ is still injective.  By the inductive hypothesis, it is also surjective.

Now let $q$ be a path in $F$ not on top at $e_1$.  Then the restriction of $q$ to $S\setm c$ is on top at some edge $e$.  By the above computation, $q$ is still on top at $e$.  Hence, the map $e\mapsto\max_eF$ is surjective.

(\ref{thm_pure_thin_2}): This statement follows from part (1) by a dual argument.

(\ref{thm_pure_thin_3}): There exists a path $r$ in $F-\{p\}$ on top at two vertical edges, say $e_1$ and $e_2$.  Let $v_1$ be the South vertex of $e_1$ and $v_2$ be the North vertex of $e_2$.  Let $e_1^{\pr}$ be the horizontal edge West of $v_1$ and $e_2^{\pr}$ the horizontal edge East of $v_2$.

Let $r^{\pr}$ be the bottom path at $e_1^{\pr}$ in $F-\{p\}$.  We claim that $r^{\pr}[v_1,v_2]=r[v_1,v_2]$ and that $r^{\pr}$ is the bottom path at $e_2^{\pr}$.

Let $v$ be the last vertex for which $r[v_1,v]=r^{\pr}[v_1,v]$.  Then $v\leq v_2$ since $r$ is the top path at $e_2$.  Let $e$ be the vertical edge South of $v$, and let $r_e$ be the top path at $e$.  Choose $v^{\pr}$ minimal such that $r_e[v^{\pr},v]=r[v^{\pr},v]$.  Since $r_e$ and $r^{\pr}$ are non-kissing $v^{\pr}\leq v_1$.  However, as $r$ is the top path at $e_1$, we either have $v_1=v^{\pr}$ or $r_e=r$.  If $r_e\neq r$, then $r_e\prec_{e_1^{\pr}}r^{\pr}$, a contradiction.  Hence $r_e=r$ and $e=e_2$.

Let $r_{e_2^{\pr}}$ be the bottom path at $e_2^{\pr}$.  Let $v$ be the smallest vertex such that $r_{e_2^{\pr}}[v,v_2]=r[v,v_2]$.  Since $r$ is the top path at $e_1$, $v_1\leq v$ holds.  If $v_1<v$, then $r_{e_2^{\pr}}$ enters $v$ from the North while $r^{\pr}$ enters from the West. However this would force $r$ and $r_{e_2^{\pr}}$ to be kissing, a contradiction.  Hence $v=v_1$ and $r^{\pr}=r_{e_2^{\pr}}$.  This completes the proof of the claim.

Define paths $q_{e_1},q_{e_2}$ such that $q_{e_1}[\cdot,v_2]=r[\cdot,v_2],\ q_{e_1}[v_2,\cdot]=r^{\pr}[v_2,\cdot]$ and $q_{e_2}[\cdot,v_2]=r^{\pr}[\cdot,v_2],\ q_{e_2}[v_2,\cdot]=r[v_2,\cdot]$.  It is easy to check that $F-\{p\}\cup\{q_{e_i}\}$ is non-kissing for $i=1,2$.  Moreover, $q_{e_1}$ and $q_{e_2}$ kiss along the unique segment $r[v_1,v_2]$.  It remains to prove that these are the only two paths that are non-kissing with $F-\{p\}$.

Let $q$ be a path such that $F-\{p\}\cup\{q\}$ is non-kissing.  Then either $q$ is on top at $e_1$ or $e_2$.

Assume $q$ is on top at $e_1$.  Let $v$ be the largest vertex for which $q[v_1,v]=r[v_1,v]$.  Since $r$ is on top at $e_2$, $v\leq v_2$ holds.  As $q$ and $r^{\pr}$ are non-kissing, we must have $v=v_2$.  If $q\neq q_{e_1}$, then they must kiss along some segment $s$.  Since $q$ is non-kissing with both $r$ and $r^{\pr}$, this segment $s$ must contain $r[v_1,v_2]$.  Since $r\prec_{e_1}q$, $q$ must enter $s$ from the West and exit South.

Let $v,v^{\pr}$ be vertices such that $s=q[v,v^{\pr}]$.  Let $e$ be the edge North of $v$, and let $p_e$ be the top path at $e$ in $F-\{p\}\cup\{q_{e_1}\}$.  Since $p_e$ and $q$ are non-kissing, we must have $p_e[v,v_1]=q[v,v_1]$.  Hence, $p_e=q_{e_1}$, a contradiction.
\end{poof}

A simplicial complex is \emph{pure} if its facets all have the same dimension.  A pure complex is \emph{thin} if every face of codimension 1 is contained in exactly two facets.  From Theorem \ref{thm_pure_thin}, we deduce the following corollary.

\begin{corollary}\label{cor_pure_thin}
For any shape $\lambda$, the reduced non-kissing complex $\wtil{\Delta}^{NK}(\lambda)$ is a pure, thin, flag simplicial complex.
\end{corollary}

This result was proven in \cite{petersen.pylyavskyy.speyer:noncrossing} and \cite{santos.stump.welker:noncrossing} for some specific shapes by identifying $\Delta^{NK}(\lambda)$ as a regular, unimodular triangulation of an order polytope with ``enough'' cone points.  The regularity of this triangulation implies that $\wtil{\Delta}^{NK}(\lambda)$ is the boundary complex of a polytope.  When $\lambda$ is a rectangle shape, this polytope is called the \emph{Grassmann Associahedron} since the triangulation reflects many of the algebraic properties of the coordinate ring of the Grassmannian, and it reduces to the usual associahedron if $\lambda$ has two rows \cite{santos.stump.welker:noncrossing}.  In the next section, we give another proof of Corollary \ref{cor_pure_thin} and of polytopality for any shape $\lambda$ by constructing $\wtil{\Delta}^{NK}(\lambda)$ from the empty complex by a sequence of suspensions and edge-stellations.

\begin{example}\label{ex_facet}
We illustrate Theorem \ref{thm_pure_thin} with the facet $F=\{145,146,236,245\}$ of $\wtil{\Delta}^{NC}_{3,6}$.  The sets in $F$ correspond to the four non-kissing paths drawn in Figure \ref{fig_kissing}.  Including the two vertical paths $234$ and $345$, each of the six paths in $F\cup\{234,345\}$ is the top path at a unique interior vertical edge.

The unique facet distinct from $F$ containing $F-\{145\}$ is $(F-\{145\})\cup\{246\}$.  If one removes $145$ from $F$, then $245$ is on top at two different vertical edges.  The segment supported by $245$ between these two vertical edges is the unique segment along which the paths $145$ and $246$ kiss.
\end{example}

The \emph{dual graph} of a pure thin complex is the set of facets where two facets are adjacent if they intersect at a codimension 1 face.  We define an orientation on the dual graph of $\wtil{\Delta}^{NK}(\lambda)$ as follows.  Let $F_1,F_2$ be adjacent facets, and let $p_1\in F_1-F_2,\ p_2\in F_2-F_1$.  Then $p_1$ and $p_2$ are kissing at a unique segment, say $p_1[v,v^{\pr}]$.  Orient the edge $F_1\ra F_2$ if $p_1$ enters $v$ from the West (equivalently, $p_1$ leaves $v^{\pr}$ to the South).  Let $\GT(\lambda)$ be the transitive closure of this relation.

\begin{theorem}[see \cite{santos.stump.welker:noncrossing}, Theorem 2.17]\label{thm_GT}
$\GT(\lambda)$ is a partially ordered set.
\end{theorem}

We call $\GT(\lambda)$ the \emph{Grid-Tamari order}.  When $\lambda$ is a $2\times n$ rectangle, $\GT(\lambda)$ is the usual Tamari lattice.  For general $\lambda$, Theorem \ref{thm_GT} is far from obvious.  In \cite{santos.stump.welker:noncrossing}, it is proved for all rectangle shapes by identifying $\GT(\lambda)$ with a poset of facets of a regular triangulation of a polytope, whose order is induced by a generic linear functional.  We establish Theorem \ref{thm_GT} as a consequence of Theorem \ref{thm_main}.

\section{Polytopal realization of the non-kissing complex}\label{sec_polytope}

If $F$ is a face of a simplicial complex $\Gamma$, the \emph{stellation} of $\Gamma$ at $F$, denoted $\st_F(\Gamma)$, is the simplicial complex $\st_F(\Gamma)=(\Gamma-F)\cup(\lk F*\pa F*\{v\})$ where $v$ is a new vertex not in the ground set of $\Gamma$.  When $\Gamma$ is the boundary complex of a simplicial polytope $P$, the stellation at $F$ may be geometrically realized by adding a new vertex to $P$ ``close'' to the center of $F$.

There are several known constructions of the dual associahedron by iteratively stellating faces of a simplex or cross-polytope \cite{lee:associahedron}, \cite{buchstaber.volodin:combinatorial}, \cite{gorsky:subword}.  We produce a similar construction for the non-kissing complex.

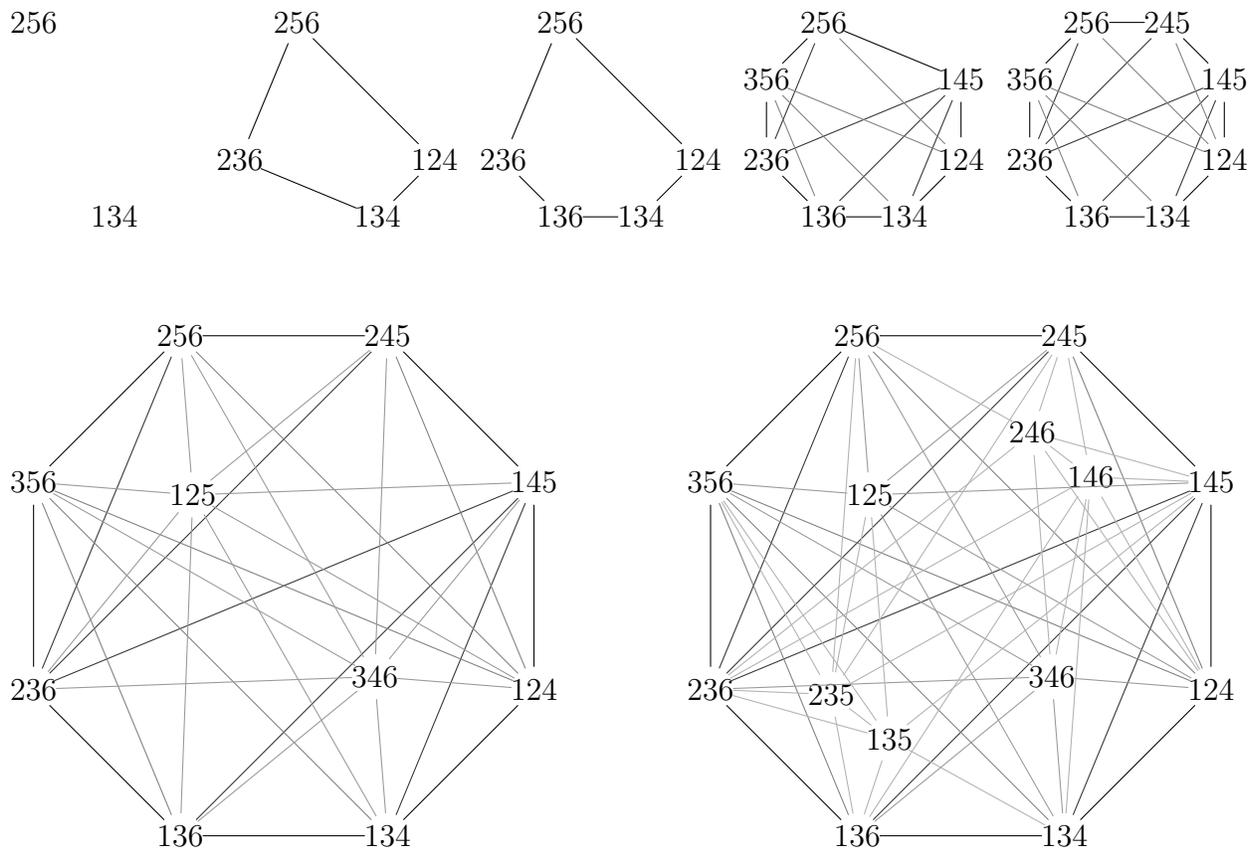
\begin{figure}
\begin{centering}
\begin{tikzpicture}[scale=.7]

\begin{scope}

\coordinate (v256) at (112.5:2cm);
\coordinate (v134) at (292.5:2cm);

\draw (v256) node{256}
      (v134) node{134};

\end{scope}

\begin{scope}[xshift=5cm]

\coordinate (v256) at (112.5:2cm);
\coordinate (v236) at (202.5:2cm);
\coordinate (v134) at (292.5:2cm);
\coordinate (v124) at (337.5:2cm);

\draw (v256) node{256}
      (v236) node{236}
      (v134) node{134}
      (v124) node{124};

\draw[shorten <=.3cm, shorten >=.3cm] (v236) -- (v134);
\draw[shorten <=.3cm, shorten >=.3cm] (v134) -- (v124);
\draw[shorten <=.3cm, shorten >=.3cm] (v236) -- (v256);
\draw[shorten <=.3cm, shorten >=.3cm] (v124) -- (v256);

\end{scope}

\begin{scope}[xshift=10cm]

\coordinate (v256) at (112.5:2cm);
\coordinate (v236) at (202.5:2cm);
\coordinate (v136) at (247.5:2cm);
\coordinate (v134) at (292.5:2cm);
\coordinate (v124) at (337.5:2cm);

\draw (v256) node{256}
      (v236) node{236}
      (v136) node{136}
      (v134) node{134}
      (v124) node{124};

\draw[shorten <=.3cm, shorten >=.3cm] (v236) -- (v136);
\draw[shorten <=.3cm, shorten >=.3cm] (v136) -- (v134);
\draw[shorten <=.3cm, shorten >=.3cm] (v134) -- (v124);
\draw[shorten <=.3cm, shorten >=.3cm] (v236) -- (v256);
\draw[shorten <=.3cm, shorten >=.3cm] (v124) -- (v256);

\end{scope}

\begin{scope}[xshift=15cm]

\coordinate (v145) at (22.5:2cm);
\coordinate (v256) at (112.5:2cm);
\coordinate (v356) at (157.5:2cm);
\coordinate (v236) at (202.5:2cm);
\coordinate (v136) at (247.5:2cm);
\coordinate (v134) at (292.5:2cm);
\coordinate (v124) at (337.5:2cm);

\draw (v256) node{256}
      (v236) node{236}
      (v136) node{136}
      (v356) node{356}
      (v134) node{134}
      (v124) node{124}
      (v145) node{145};

\draw[shorten <=.3cm, shorten >=.3cm] (v145) -- (v256);
\draw[shorten <=.3cm, shorten >=.3cm] (v256) -- (v356);
\draw[shorten <=.3cm, shorten >=.3cm] (v356) -- (v236);
\draw[shorten <=.3cm, shorten >=.3cm] (v236) -- (v136);
\draw[shorten <=.3cm, shorten >=.3cm] (v136) -- (v134);
\draw[shorten <=.3cm, shorten >=.3cm] (v134) -- (v124);
\draw[shorten <=.3cm, shorten >=.3cm] (v124) -- (v145);

\draw[black!80, shorten <=.3cm, shorten >=.3cm] (v145) -- (v236);
\draw[black!80, shorten <=.3cm, shorten >=.3cm] (v145) -- (v136);
\draw[black!80, shorten <=.3cm, shorten >=.3cm] (v145) -- (v134);
\draw[black!80, shorten <=.3cm, shorten >=.3cm] (v236) -- (v256);
\draw[black!50, shorten <=.3cm, shorten >=.3cm] (v356) -- (v136);
\draw[black!50, shorten <=.3cm, shorten >=.3cm] (v356) -- (v134);
\draw[black!50, shorten <=.3cm, shorten >=.3cm] (v356) -- (v124);
\draw[black!50, shorten <=.3cm, shorten >=.3cm] (v124) -- (v256);

\end{scope}

\begin{scope}[xshift=20cm]

\coordinate (v145) at (22.5:2cm);
\coordinate (v245) at (67.5:2cm);
\coordinate (v256) at (112.5:2cm);
\coordinate (v356) at (157.5:2cm);
\coordinate (v236) at (202.5:2cm);
\coordinate (v136) at (247.5:2cm);
\coordinate (v134) at (292.5:2cm);
\coordinate (v124) at (337.5:2cm);

\draw (v256) node{256}
      (v236) node{236}
      (v136) node{136}
      (v356) node{356}
      (v134) node{134}
      (v124) node{124}
      (v145) node{145}
      (v245) node{245};

\draw[shorten <=.3cm, shorten >=.3cm] (v145) -- (v245);
\draw[shorten <=.3cm, shorten >=.3cm] (v245) -- (v256);
\draw[shorten <=.3cm, shorten >=.3cm] (v256) -- (v356);
\draw[shorten <=.3cm, shorten >=.3cm] (v356) -- (v236);
\draw[shorten <=.3cm, shorten >=.3cm] (v236) -- (v136);
\draw[shorten <=.3cm, shorten >=.3cm] (v136) -- (v134);
\draw[shorten <=.3cm, shorten >=.3cm] (v134) -- (v124);
\draw[shorten <=.3cm, shorten >=.3cm] (v124) -- (v145);

\draw[black!80, shorten <=.3cm, shorten >=.3cm] (v145) -- (v236);
\draw[black!80, shorten <=.3cm, shorten >=.3cm] (v145) -- (v136);
\draw[black!80, shorten <=.3cm, shorten >=.3cm] (v145) -- (v134);
\draw[black!80, shorten <=.3cm, shorten >=.3cm] (v236) -- (v245);
\draw[black!80, shorten <=.3cm, shorten >=.3cm] (v236) -- (v256);
\draw[black!50, shorten <=.3cm, shorten >=.3cm] (v356) -- (v136);
\draw[black!50, shorten <=.3cm, shorten >=.3cm] (v356) -- (v134);
\draw[black!50, shorten <=.3cm, shorten >=.3cm] (v356) -- (v124);
\draw[black!50, shorten <=.3cm, shorten >=.3cm] (v124) -- (v256);
\draw[black!50, shorten <=.3cm, shorten >=.3cm] (v124) -- (v245);

\end{scope}

\end{tikzpicture}

\vspace{1cm}

\begin{tikzpicture}[scale=1.8]

\begin{scope}

\coordinate (v145) at (22.5:2cm);
\coordinate (v245) at (67.5:2cm);
\coordinate (v256) at (112.5:2cm);
\coordinate (v356) at (157.5:2cm);
\coordinate (v236) at (202.5:2cm);
\coordinate (v136) at (247.5:2cm);
\coordinate (v134) at (292.5:2cm);
\coordinate (v124) at (337.5:2cm);

\coordinate (v346) at (315:.95cm);
\coordinate (v125) at (135:.95cm);

\draw (v256) node{256}
      (v236) node{236}
      (v136) node{136}
      (v356) node{356}
      (v134) node{134}
      (v124) node{124}
      (v145) node{145}
      (v245) node{245}
      (v125) node{125}
      (v346) node{346};

\draw[shorten <=.3cm, shorten >=.3cm] (v145) -- (v245);
\draw[shorten <=.3cm, shorten >=.3cm] (v245) -- (v256);
\draw[shorten <=.3cm, shorten >=.3cm] (v256) -- (v356);
\draw[shorten <=.3cm, shorten >=.3cm] (v356) -- (v236);
\draw[shorten <=.3cm, shorten >=.3cm] (v236) -- (v136);
\draw[shorten <=.3cm, shorten >=.3cm] (v136) -- (v134);
\draw[shorten <=.3cm, shorten >=.3cm] (v134) -- (v124);
\draw[shorten <=.3cm, shorten >=.3cm] (v124) -- (v145);

\draw[black!80, shorten <=.3cm, shorten >=.3cm] (v145) -- (v236);
\draw[black!80, shorten <=.3cm, shorten >=.3cm] (v145) -- (v136);
\draw[black!80, shorten <=.3cm, shorten >=.3cm] (v145) -- (v134);
\draw[black!80, shorten <=.3cm, shorten >=.3cm] (v236) -- (v245);
\draw[black!80, shorten <=.3cm, shorten >=.3cm] (v236) -- (v256);
\draw[black!50, shorten <=.3cm, shorten >=.3cm] (v356) -- (v136);
\draw[black!50, shorten <=.3cm, shorten >=.3cm] (v356) -- (v134);
\draw[black!50, shorten <=.3cm, shorten >=.3cm] (v356) -- (v124);
\draw[black!50, shorten <=.3cm, shorten >=.3cm] (v124) -- (v256);
\draw[black!50, shorten <=.3cm, shorten >=.3cm] (v124) -- (v245);

\draw[black!40, shorten <=.3cm, shorten >=.3cm] (v125) -- (v134);
\draw[black!40, shorten <=.3cm, shorten >=.3cm] (v125) -- (v124);
\draw[black!40, shorten <=.3cm, shorten >=.3cm] (v125) -- (v145);
\draw[black!40, shorten <=.3cm, shorten >=.3cm] (v125) -- (v245);
\draw[black!40, shorten <=.3cm, shorten >=.3cm] (v125) -- (v256);
\draw[black!40, shorten <=.3cm, shorten >=.3cm] (v125) -- (v356);
\draw[black!40, shorten <=.3cm, shorten >=.3cm] (v125) -- (v236);
\draw[black!40, shorten <=.3cm, shorten >=.3cm] (v125) -- (v136);

\draw[black!40, shorten <=.3cm, shorten >=.3cm] (v346) -- (v256);
\draw[black!40, shorten <=.3cm, shorten >=.3cm] (v346) -- (v356);
\draw[black!40, shorten <=.3cm, shorten >=.3cm] (v346) -- (v236);
\draw[black!40, shorten <=.3cm, shorten >=.3cm] (v346) -- (v136);
\draw[black!40, shorten <=.3cm, shorten >=.3cm] (v346) -- (v134);
\draw[black!40, shorten <=.3cm, shorten >=.3cm] (v346) -- (v124);
\draw[black!40, shorten <=.3cm, shorten >=.3cm] (v346) -- (v145);
\draw[black!40, shorten <=.3cm, shorten >=.3cm] (v346) -- (v245);

\end{scope}

\begin{scope}[xshift=5cm]

\coordinate (v145) at (22.5:2cm);
\coordinate (v245) at (67.5:2cm);
\coordinate (v256) at (112.5:2cm);
\coordinate (v356) at (157.5:2cm);
\coordinate (v236) at (202.5:2cm);
\coordinate (v136) at (247.5:2cm);
\coordinate (v134) at (292.5:2cm);
\coordinate (v124) at (337.5:2cm);

\coordinate (v346) at (315:.95cm);
\coordinate (v125) at (135:.95cm);

\coordinate (v135) at (245:1.25cm);
\coordinate (v235) at (220:1.25cm);
\coordinate (v146) at (40:1.25cm);
\coordinate (v246) at (65:1.25cm);

\draw (v256) node{256}
      (v236) node{236}
      (v136) node{136}
      (v356) node{356}
      (v134) node{134}
      (v124) node{124}
      (v145) node{145}
      (v245) node{245}
      (v125) node{125}
      (v346) node{346}
      (v135) node{135}
      (v235) node{235}
      (v146) node{146}
      (v246) node{246};

\draw[shorten <=.3cm, shorten >=.3cm] (v145) -- (v245);
\draw[shorten <=.3cm, shorten >=.3cm] (v245) -- (v256);
\draw[shorten <=.3cm, shorten >=.3cm] (v256) -- (v356);
\draw[shorten <=.3cm, shorten >=.3cm] (v356) -- (v236);
\draw[shorten <=.3cm, shorten >=.3cm] (v236) -- (v136);
\draw[shorten <=.3cm, shorten >=.3cm] (v136) -- (v134);
\draw[shorten <=.3cm, shorten >=.3cm] (v134) -- (v124);
\draw[shorten <=.3cm, shorten >=.3cm] (v124) -- (v145);

\draw[black!80, shorten <=.3cm, shorten >=.3cm] (v145) -- (v236);
\draw[black!80, shorten <=.3cm, shorten >=.3cm] (v145) -- (v136);
\draw[black!80, shorten <=.3cm, shorten >=.3cm] (v145) -- (v134);
\draw[black!80, shorten <=.3cm, shorten >=.3cm] (v236) -- (v245);
\draw[black!80, shorten <=.3cm, shorten >=.3cm] (v236) -- (v256);
\draw[black!50, shorten <=.3cm, shorten >=.3cm] (v356) -- (v136);
\draw[black!50, shorten <=.3cm, shorten >=.3cm] (v356) -- (v134);
\draw[black!50, shorten <=.3cm, shorten >=.3cm] (v356) -- (v124);
\draw[black!50, shorten <=.3cm, shorten >=.3cm] (v124) -- (v256);
\draw[black!50, shorten <=.3cm, shorten >=.3cm] (v124) -- (v245);

\draw[black!40, shorten <=.3cm, shorten >=.3cm] (v125) -- (v134);
\draw[black!40, shorten <=.3cm, shorten >=.3cm] (v125) -- (v124);
\draw[black!40, shorten <=.3cm, shorten >=.3cm] (v125) -- (v145);
\draw[black!40, shorten <=.3cm, shorten >=.3cm] (v125) -- (v245);
\draw[black!40, shorten <=.3cm, shorten >=.3cm] (v125) -- (v256);
\draw[black!40, shorten <=.3cm, shorten >=.3cm] (v125) -- (v356);
\draw[black!40, shorten <=.3cm, shorten >=.3cm] (v346) -- (v256);
\draw[black!40, shorten <=.3cm, shorten >=.3cm] (v346) -- (v356);
\draw[black!40, shorten <=.3cm, shorten >=.3cm] (v346) -- (v236);
\draw[black!40, shorten <=.3cm, shorten >=.3cm] (v346) -- (v136);
\draw[black!40, shorten <=.3cm, shorten >=.3cm] (v346) -- (v134);
\draw[black!40, shorten <=.3cm, shorten >=.3cm] (v346) -- (v124);

\draw[black!30, shorten <=.3cm, shorten >=.3cm] (v235) -- (v125);
\draw[black!30, shorten <=.3cm, shorten >=.3cm] (v235) -- (v145);
\draw[black!30, shorten <=.3cm, shorten >=.3cm] (v235) -- (v245);
\draw[black!30, shorten <=.3cm, shorten >=.3cm] (v235) -- (v256);
\draw[black!30, shorten <=.3cm, shorten >=.3cm] (v235) -- (v356);
\draw[black!30, shorten <=.3cm, shorten >=.3cm] (v235) -- (v236);
\draw[black!30, shorten <=.3cm, shorten >=.3cm] (v235) -- (v136);

\draw[black!30, shorten <=.3cm, shorten >=.3cm] (v135) -- (v125);
\draw[black!30, shorten <=.3cm, shorten >=.3cm] (v135) -- (v136);
\draw[black!30, shorten <=.3cm, shorten >=.3cm] (v135) -- (v134);
\draw[black!30, shorten <=.3cm, shorten >=.3cm] (v135) -- (v145);
\draw[black!30, shorten <=.3cm, shorten >=.3cm] (v135) -- (v356);
\draw[black!30, shorten <=.3cm, shorten >=.3cm] (v135) -- (v236);
\draw[black!30, shorten <=.3cm, shorten >=.3cm] (v135) -- (v235);

\draw[black!30, shorten <=.3cm, shorten >=.3cm] (v146) -- (v346);
\draw[black!30, shorten <=.3cm, shorten >=.3cm] (v146) -- (v236);
\draw[black!30, shorten <=.3cm, shorten >=.3cm] (v146) -- (v136);
\draw[black!30, shorten <=.3cm, shorten >=.3cm] (v146) -- (v134);
\draw[black!30, shorten <=.3cm, shorten >=.3cm] (v146) -- (v124);
\draw[black!30, shorten <=.3cm, shorten >=.3cm] (v146) -- (v145);
\draw[black!30, shorten <=.3cm, shorten >=.3cm] (v146) -- (v245);

\draw[black!30, shorten <=.3cm, shorten >=.3cm] (v246) -- (v346);
\draw[black!30, shorten <=.3cm, shorten >=.3cm] (v246) -- (v245);
\draw[black!30, shorten <=.3cm, shorten >=.3cm] (v246) -- (v256);
\draw[black!30, shorten <=.3cm, shorten >=.3cm] (v246) -- (v236);
\draw[black!30, shorten <=.3cm, shorten >=.3cm] (v246) -- (v124);
\draw[black!30, shorten <=.3cm, shorten >=.3cm] (v246) -- (v145);
\draw[black!30, shorten <=.3cm, shorten >=.3cm] (v246) -- (v146);

\end{scope}

\end{tikzpicture}
\caption{\scriptsize\label{fig_stellation}A construction of the reduced non-crossing complex $\wtil{\Delta}^{NC}_{3,6}$ by a sequence of suspensions and edge stellations.}
\end{centering}
\end{figure}

Fix a shape $\lambda$, and let $c$ be a SE-corner of $\lambda$.  Let $v$ be the point one step NW of $c$.  If $v$ is not an interior vertex of $\lambda$, then $\wtil{\Delta}^{NK}(\lambda\setm c)=\wtil{\Delta}^{NK}(\lambda)$.  On the other hand, if $v$ is interior, we construct a sequence of complexes $\Gamma_0,\ldots,\Gamma_l$ such that
\begin{itemize}
\item $\Gamma_0$ is isomorphic to the suspension of $\wtil{\Delta}^{NK}(\lambda\setm c)$,
\item $\Gamma_l=\wtil{\Delta}^{NK}(\lambda)$, and
\item $\Gamma_i$ is the stellation of $\Gamma_{i-1}$ at some edge for all $i$.
\end{itemize}

Given a path in $\lambda\setm c$, we extend it (uniquely) to a path in $\lambda$ that does not turn at $v$.  Then two paths in $\lambda\setm c$ are non-kissing if and only if their extensions to $\lambda$ are non-kissing.  Let $\Gamma_0$ be the suspension of $\wtil{\Delta}^{NK}(\lambda\setm c)$ where the two new vertices correspond to the two paths $q_W,q_N$ that only turn at $v$, where $q_W$ enters $v$ from the West and $q_N$ enters $v$ from the North.

Let $e_W$ be the horizontal edge West of $v$, and let $e_N$ be the vertical edge North of $v$.  Let $p_1,\ldots,p_k$ be the list of paths distinct from $q_W$ that turn at $v$ and contain $e_W$, ordered so that if $i<j\leq k$ then $p_i\prec_{e_W}p_j$.  This is well-defined since $<_{e_W}$ is a total order on these paths.  Similarly, let $p_{k+1},\ldots,p_l$ be the list of paths that turn at $v$ and contain $e_N$, ordered so that if $k<i<j$ then $p_i\succ_{e_N}p_j$.  For each $i$, let $r_i$ be the same path as $p_i$ except that it continues straight through $v$.

Then for each $i\leq k$, define $\Gamma_i$ recursively as the complex $\st_{\{r_i,q_W\}}(\Gamma_{i-1})$, where the new vertex is labeled $p_i$.  For $i>k$, we define $\Gamma_i$ as the complex $\st_{\{r_i,q_N\}}(\Gamma_{i-1})$, where the new vertex is again labeled $p_i$.

With the above set-up, the following result is elementary, if somewhat tedious to verify.

\begin{theorem}\label{thm_non-kissing_stellation}
$\Gamma_l=\wtil{\Delta}^{NK}(\lambda)$.
\end{theorem}

\begin{poof}
Let $p,q$ be two paths supported by $\lambda$.  We prove that $p$ and $q$ are adjacent in $\Gamma_l$ if and only if they are non-kissing.

If neither $p$ nor $q$ turns at $v$, then $p$ and $q$ are adjacent in $\Gamma_0$ if and only if they are non-kissing.  As these edges are not stellated by the construction, it follows that $p$ and $q$ are adjacent in $\Gamma_l$ exactly when they are non-kissing.

Assume $q=q_W$.  Then $q$ kisses $p$ only if $p$ leaves $v$ to the East.  If $p$ and $q$ kiss at $(v)$, then either $p=q_N$ or $p=q_i$ for some $i>k$.  In either case, they are not adjacent in $\Gamma_l$.  If $p$ and $q$ kiss at a segment $s$ containing $e_W$, then $p=r_i$ for some $i\leq k$.  In this case, $p$ and $q$ are separated in $\Gamma_i$.  As $q$ is adjacent to every other vertex of $\Gamma_l$, we are done in this case.  A similar argument holds if $q=q_N$.

Now assume $p=p_i$ for some $i\leq k$.  Suppose $p$ and $q$ kiss along a segment $s$ not containing $v$.  Then $r_i$ and $q$ also kiss along $s$.  If $q$ does not turn at $v$, then $r_i$ and $q$ are not adjacent in $\Gamma_0$, so $p$ and $q$ are not adjacent in $\Gamma_i$.  If $q$ does turn at $v$, then $q=p_j$.  Without loss of generality, we may assume $i<j$.  Then $r_i$ and $r_j$ are not adjacent, so $p$ and $r_j$ are not adjacent in $\Gamma_i$ and $p$ and $q$ are not adjacent in $\Gamma_j$.

Assume $p=p_i$ for some $i\leq k$ and suppose $p$ and $q$ only kiss along a segment $s$ containing $v$.  If $s=(v)$, then $q=p_j$ for some $j>k$ or $q=q_N$.  In either case, $p$ and $q$ are not adjacent in $\Gamma_l$.  If $s$ contains $e_W$ then $q=r_j$ for some $j\leq k$.  As $p_j\prec_{e_W}p_i$ we deduce that $j<i$.  Hence, $q_W$ is not adjacent to $q$ in $\Gamma_{i-1}$, so $p$ and $q$ are not adjacent in $\Gamma_i$.

Now assume $p=p_i$ for some $i\leq k$ and suppose $p$ and $q$ are non-kissing.  If $q$ does not contain $v$, then $q$ is adjacent to $r_i$ and $q_W$ in $\Gamma_0$, so $p$ and $q$ are adjacent in $\Gamma_i$.  If $q$ contains $v$, then either $q=q_W$, $q=r_j$ for some $j>k$, $q=r_j$ for some $j\leq k$, or $q=p_j$ for some $j\leq k$.  The first case has already been handled.  In the second case, $r_j$ and $r_i$ are non-kissing, so $p$ and $q$ are adjacent in $\Gamma_i$.  In the third case, either $r_i$ and $r_j$ are non-kissing, or $i<j$; for both situations, $p$ and $q$ are adjacent in $\Gamma_i$.  Finally, if $q=p_j$ for some $j\leq k$, we may assume $i<j$ without loss of generality.  Then $p$ and $r_j$ are adjacent in $\Gamma_i$, so $p$ and $q$ are adjacent in $\Gamma_j$.

A similar argument holds when $p=p_i$ and $i>k$.  This completes the proof.
\end{poof}

\section{Lattice properties of biclosed sets}\label{sec_biclosed}

A \emph{closure operator} on a set $S$ is an operator $X\mapsto\ov{X}$ on subsets of $S$ such that for $X,Y\subseteq S$,
\begin{align*}
& X\subseteq\ov{X},\\
& \ov{\ov{X}}=\ov{X},\ \mbox{and}\\
& X\subseteq Y\ \mbox{implies}\ \ov{X}\subseteq\ov{Y}.
\end{align*}

In addition, we assume $\ov{\emptyset}=\emptyset$.  A subset $X$ of $S$ is \emph{closed} if $X=\ov{X}$.  A set $X$ is \emph{co-closed} (or \emph{open}) if $S-X$ is closed.  We say $X$ is \emph{biclosed} if $X$ and $S-X$ are both closed.  We let $\Bic(S)$ be the poset of biclosed subsets of $S$ ordered by inclusion.  By our assumption, $S$ and $\emptyset$ are always biclosed.

Two important families of closure operators are the convex closure and rank-2 convex closure on a finite subset of $\Rbb^n$.  Given a finite subset $S$ of $\Rbb^n$, the \emph{convex closure} of a subset $X$ of $S$ is the set of points in $S$ that can be expressed as convex linear combinations of points in $X$.  The \emph{rank-2 convex closure} (or \emph{2-closure}) of $X$ is the smallest subset $\ov{X}$ of $S$ containing $X$ such that if $x,y\in\ov{X}$ and $z\in S$ such that $z=\lambda x+(1-\lambda)y$ for some $\lambda\in[0,1]$ then $z\in\ov{X}$; that is, $\ov{X}$ is convex along lines.

\begin{remark}
What we call biclosed sets are often called clopen sets elsewhere in the literature; see, for example \cite{santocanale.wehrung:lattices}.  The term biclosed typically refers to a subset of a convex geometry which is 2-closed and whose complement is 2-closed.  We choose the term biclosed because all of the closure operators we consider come from some convex geometry in this way.
\end{remark}

\subsection{Semidistributive lattices}

A collection $\Bcal$ of subsets of $S$ is \emph{ordered by single-step inclusion} if for all $X,Y\in\Bcal$ such that $X\subsetneq Y$ there exists $y\in Y-X$ such that $X\cup\{y\}\in\Bcal$.  If $\emptyset,S\in\Bcal$ and $\Bcal$ is ordered by single-step inclusion, then it is a graded lattice with rank function $X\mapsto|X|$ for $X\in\Bcal$; in particular, every maximal chain has length $|S|$.

A lattice $L$ is \emph{meet-semidistributive} if $L$ satisfies $x\wedge z=y\wedge z\ \Ra\ (x\vee y)\wedge z=x\wedge z$ for $x,y,z\in L$.  A lattice is \emph{join-semidistributive} if its dual is meet-semidistributive.  A lattice is \emph{semidistributive} if it is both meet- and join-semidistributive.

For fixed $z\in L$, the map $x\mapsto x\wedge z$ is an order-preserving map $L\ra L$ that preserves meets.  Thus, Lemma \ref{lem_local_map_join}(\ref{lem_lmj_2}) determines a local test for meet-semidistributivity.

\begin{theorem}\label{thm_semidistributive}
Let $S$ be a set with a closure operator.  If
\begin{enumerate}
\item\label{thm_semi_1} $\Bic(S)$ is ordered by single-step inclusion, and
\item\label{thm_semi_2} $W\cup\ov{(X\cup Y)-W}$ is biclosed for $W,X,Y\in\Bic(S)$ with $W\subseteq X\cap Y$,
\end{enumerate}
then $\Bic(S)$ is a semidistributive lattice.
\end{theorem}

\begin{poof}
If $W,X,Y\in\Bic(S)$ with $W\subseteq X\cap Y$, then
$$X\cup Y\subseteq W\cup\ov{(X\cup Y)-W}\subseteq\ov{X\cup Y},$$
so $X\vee Y$ and $W\cup\ov{(X\cup Y)-W}$ are equal if the latter is biclosed.  Taking $W=\emptyset$, condition (\ref{thm_semi_2}) implies $\Bic(S)$ is a lattice.

Since $\Bic(S)$ is a self-dual poset, semidistributivity follows from meet-semidistributivity.  By the above discussion, it suffices to show for $W,X,Y,Z\in\Bic(S)$ if $X$ and $Y$ both cover $W$ and $X\wedge Z=Y\wedge Z$, then $(X\vee Y)\wedge Z=X\wedge Z$.

By (\ref{thm_semi_1}), there exists $s,t\in S$ such that $X=W\cup\{s\}$ and $Y=W\cup\{t\}$.  By (\ref{thm_semi_2}), $X\vee Y=W\cup\ov{\{s,t\}}$.  If $W\wedge Z<(X\vee Y)\wedge Z$, then there exists $u\in(X\vee Y)\wedge Z$ such that $(W\wedge Z)\cup\{u\}$ is biclosed.  Then $u$ is an element of $(X\vee Y)-W$, so $u\in\ov{\{s,t\}}$.  Since $W\wedge Z=X\wedge Z=Y\wedge Z$, the elements $s,t$ are not in $W\wedge Z$ and $u\neq s,\ u\neq t$.  However, this implies $\{s,t\}$ is contained in the complement of $(W\wedge Z)\cup\{u\}$, contradicting the assumption that this set is biclosed.  Hence, $W\wedge Z=(X\vee Y)\wedge Z$ holds.
\end{poof}

\begin{example}\label{ex_weak}
The weak order on permutations may be identified with a collection of ``biclosed'' subsets of $\binom{[n]}{2}$, ordered by inclusion.  A subset $X$ of $\binom{[n]}{2}$ is closed if $\{i,k\}$ is in $X$ whenever $\{i,j\}$ and $\{j,k\}$ are in $X$ for some $j$ with $i<j<k$.  Then $X$ is biclosed if both $X$ and $\binom{[n]}{2}-X$ are closed.  The map taking a permutation to its inversion set is an isomorphism between the weak order and the poset of biclosed subsets of $\binom{[n]}{2}$.

More generally, the weak order on any finite Coxeter group may be identified with a poset of biclosed sets of positive roots ordered by inclusion.  That these posets are ordered by single-step inclusion is well-known.  Dyer proved that $W\cup\ov{(X\cup Y)-W}$ is a biclosed set whenever $W,X,Y$ are biclosed and $W\subseteq X\cap Y$ \cite{dyer:weak}.  He also proved this holds for infinite root systems if $\ov{X\cup Y}$ is finite.  By Theorem \ref{thm_semidistributive} we may deduce that the weak order for finite Coxeter groups is a semidistributive lattice.  Other proofs of semidistributivity appear in \cite{poly-barbut:treillis} and \cite{reading:lattice}.
\end{example}

\subsection{Congruence-normal and congruence-uniform lattices}\label{subsec_CN}

A subset $C$ of a poset $P$ is \emph{order-convex} if $z\in C$ whenever $x,y\in C$ and $x\leq z\leq y$.  Given an order-convex subset $C$ of $P$, the doubling $P[C]$ is the induced subposet of $P\times\{0,1\}$ with elements
$$P[C]=(P_{\leq C}\times\{0\})\sqcup[(P-P_{\leq C})\cup C]\times\{1\},$$
where $P_{\leq C}=\{x\in P:\ (\exists c\in C)\ x\leq c\}$.  If $P$ is a lattice, then $P[C]$ is a lattice where
$$(x,\epsilon)\vee(y,\epsilon^{\pr})=\begin{cases}(x\vee y,\max(\epsilon,\epsilon^{\pr}))\ &\mbox{if }x\vee y\in P_{\leq C}\\(x\vee y,1)\ &\mbox{otherwise}\end{cases},$$
for $(x,\epsilon),(y,\epsilon^{\pr})\in P[C]$.  A finite lattice $L$ is \emph{congruence-normal} if there exists a sequence of lattices $L_1,\ldots,L_l$ such that $L_1$ is the one-element lattice, $L_l=L$, and for all $i$, there exists an order convex subset $C_i$ of $L_i$ such that $L_{i+1}\cong L_i[C_i]$.  A lattice is \emph{congruence-uniform} if it is both congruence-normal and semidistributive.

In Section \ref{sec_lattices}, we defined congruence-uniformity in terms of lattice congruences.  The equivalence of these two definitions was proved by Day \cite{day:congruence}.  Additionally, congruence-uniform lattices may be characterized as lattice quotients of free lattices for which every fiber is a closed interval.  As free lattices are typically infinite, this interval property is quite special.

\begin{figure}
\begin{centering}
\begin{tikzpicture}[scale=.7]

\begin{scope}
\filldraw (0,0) circle(.5mm);
\draw[black!70] (0,0) circle(2mm);
\end{scope}

\begin{scope}[xshift=.8cm,yshift=-.75cm]
\draw (0,0) -- (0,1.5);
\draw[black!70] (0,.75) ellipse[x radius=.2cm, y radius=.9cm];
\end{scope}

\begin{scope}[xshift=2cm,yshift=-1cm]
\draw (0,0) -- (-.5,1) -- (0,2) -- (.5,1) -- cycle;
\draw[black!70] (-.5,1) circle(2mm);
\draw[black!70] (.5,1) circle(2mm);
\end{scope}

\begin{scope}[xshift=4cm,yshift=-1.3cm]
\draw (0,0) -- (-.6,.8) -- (-.6,1.8) -- (0,2.6) -- (.6,1.8) -- (.6,.8) -- cycle;
\draw[black!70] (0,1.3) ellipse[x radius=1cm, y radius=1.5cm];
\end{scope}

\begin{scope}[xshift=8cm,yshift=-2.1cm]
  \coordinate (O) at (0,0);
  \coordinate (A) at (-1.6,0.9);  
  \coordinate (B) at (-0.8,1.1);
  \coordinate (D) at (0.2,1.3);
  \coordinate (F) at (1.6,0.9);

  \draw (O) -- (A);
  \draw (O) -- (D);
  \draw (O) -- (F);
  \draw (A) -- ($(A)+(B)$);
  \draw (D) -- ($(B)+(D)$);
  \draw ($(A)+(B)$) -- ($(A)+(B)+(D)$);
  \draw ($(B)+(D)$) -- ($(A)+(B)+(D)$);
  \draw (D) -- ($(B)+(D)$);
  \draw (D) -- ($(D)+(F)$);
  \draw (F) -- ($(D)+(F)$);
  \draw ($(A)+(B)$) -- ($(A)+(B)+(D)$);
  \draw ($(B)+(D)$) -- ($(A)+(B)+(D)$);
  \draw ($(B)+(D)$) -- ($(B)+(D)+(F)$);
  \draw ($(D)+(F)$) -- ($(B)+(D)+(F)$);
  \draw ($(A)+(B)+(D)$) -- ($(A)+(B)+(D)+(F)$);
  \draw ($(B)+(D)+(F)$) -- ($(A)+(B)+(D)+(F)$);
  \draw (A) -- ($(A)+(F)$);
  \draw (F) -- ($(A)+(F)$);
  \draw (F) -- ($(D)+(F)$);
  \draw (F) -- ($(A)+(F)$);
  \draw ($(A)+(F)$) -- ($(A)+(B)+(F)$);
  \draw ($(A)+(B)$) -- ($(A)+(B)+(F)$);
  \draw ($(A)+(F)$) -- ($(A)+(B)+(F)$);
  \draw ($(A)+(B)+(F)$) -- ($(A)+(B)+(D)+(F)$);

  \draw[black!70] (-.3,2.1) ellipse[rotate=150,x radius=2.5cm,y radius=.5cm];
\end{scope}

\begin{scope}[xshift=13cm,yshift=-2.4cm]
  \coordinate (O) at (0,0);
  \coordinate (A) at (-1.6,0.8);  
  \coordinate (B) at (-0.8,1);
  \coordinate (D) at (0.2,1.2);
  \coordinate (E) at (0.8,1);
  \coordinate (F) at (1.6,0.8);

  \draw (O) -- (A);
  \draw (O) -- (D);
  \draw (O) -- (F);
  \draw (A) -- ($(A)+(B)$);
  \draw (D) -- ($(B)+(D)$);
  \draw (D) -- ($(D)+(E)$);
  \draw (F) -- ($(E)+(F)$);
  \draw ($(A)+(B)$) -- ($(A)+(B)+(D)$);
  \draw ($(B)+(D)$) -- ($(A)+(B)+(D)$);
  \draw ($(B)+(D)$) -- ($(B)+(D)+(E)$);
  \draw ($(D)+(E)$) -- ($(B)+(D)+(E)$);
  \draw ($(D)+(E)$) -- ($(D)+(E)+(F)$);
  \draw ($(E)+(F)$) -- ($(D)+(E)+(F)$);
  \draw ($(A)+(B)$) -- ($(A)+(B)+(D)$);
  \draw ($(A)+(B)+(D)$) -- ($(A)+(B)+(D)+(E)$);
  \draw ($(B)+(D)+(E)$) -- ($(A)+(B)+(D)+(E)$);
  \draw ($(B)+(D)+(E)$) -- ($(B)+(D)+(E)+(F)$);
  \draw ($(D)+(E)+(F)$) -- ($(B)+(D)+(E)+(F)$);
  \draw ($(A)+(B)+(D)+(E)$) -- ($(A)+(B)+(D)+(E)+(F)$);
  \draw ($(B)+(D)+(E)+(F)$) -- ($(A)+(B)+(D)+(E)+(F)$);
  \draw (A) -- ($(A)+(F)$);
  \draw (F) -- ($(A)+(F)$);
  \draw ($(E)+(F)$) -- ($(D)+(E)+(F)$);
  \draw ($(E)+(F)$) -- ($(A)+(E)+(F)$);
  \draw ($(A)+(F)$) -- ($(A)+(E)+(F)$);
  \draw ($(A)+(F)$) -- ($(A)+(B)+(F)$);
  \draw ($(A)+(B)$) -- ($(A)+(B)+(F)$);
  \draw ($(A)+(B)+(F)$) -- ($(A)+(B)+(E)+(F)$);
  \draw ($(A)+(E)+(F)$) -- ($(A)+(B)+(E)+(F)$);
  \draw ($(A)+(B)+(E)+(F)$) -- ($(A)+(B)+(D)+(E)+(F)$);

  \draw[black!70] ($(A)+(F)$) circle(2mm);
  \draw[black!70] ($(B)+(D)+(E)$) circle(2mm);
  \draw[black!70] (-2.3,2.4) ellipse[rotate=83,x radius=.8cm,y radius=.3cm];
  \draw[black!70] (2.5,2.4) ellipse[rotate=83,x radius=.8cm,y radius=.3cm];  
\end{scope}

\begin{scope}[xshift=19cm,yshift=-2.4cm]
  \coordinate (O) at (0,0);
  \coordinate (A) at (-1.6,0.6);  
  \coordinate (B) at (-0.8,0.8);
  \coordinate (C) at (-0.1,1);
  \coordinate (D) at (0.1,1);
  \coordinate (E) at (0.8,0.8);
  \coordinate (F) at (1.6,0.6);
  \draw (O) -- (A);
  \draw (O) -- (D);
  \draw (O) -- (F);
  \draw (A) -- ($(A)+(B)$);
  \draw (D) -- ($(B)+(D)$);
  \draw (D) -- ($(D)+(E)$);
  \draw (F) -- ($(E)+(F)$);
  \draw ($(A)+(B)$) -- ($(A)+(B)+(C)$);
  \draw ($(A)+(B)$) -- ($(A)+(B)+(D)$);
  \draw ($(B)+(D)$) -- ($(A)+(B)+(D)$);
  \draw ($(B)+(D)$) -- ($(B)+(D)+(E)$);
  \draw ($(D)+(E)$) -- ($(B)+(D)+(E)$);
  \draw ($(D)+(E)$) -- ($(D)+(E)+(F)$);
  \draw ($(E)+(F)$) -- ($(D)+(E)+(F)$);
  \draw ($(D)+(E)+(F)$) -- ($(C)+(D)+(E)+(F)$);
  \draw ($(B)+(D)+(E)$) -- ($(B)+(C)+(D)+(E)$);
  \draw ($(A)+(B)+(D)$) -- ($(A)+(B)+(C)+(D)$);
  \draw ($(A)+(B)+(C)$) -- ($(A)+(B)+(C)+(D)$);
  \draw ($(A)+(B)+(C)+(D)$) -- ($(A)+(B)+(C)+(D)+(E)$);
  \draw ($(B)+(C)+(D)+(E)$) -- ($(A)+(B)+(C)+(D)+(E)$);
  \draw ($(B)+(C)+(D)+(E)$) -- ($(B)+(C)+(D)+(E)+(F)$);
  \draw ($(C)+(D)+(E)+(F)$) -- ($(B)+(C)+(D)+(E)+(F)$);
  \draw ($(A)+(B)+(C)+(D)+(E)$) -- ($(A)+(B)+(C)+(D)+(E)+(F)$);
  \draw ($(B)+(C)+(D)+(E)+(F)$) -- ($(A)+(B)+(C)+(D)+(E)+(F)$);
  \draw (A) -- ($(A)+(F)$);
  \draw (F) -- ($(A)+(F)$);
  \draw ($(A)+(F)$) -- ($(A)+(C)+(F)$);
  \draw ($(E)+(F)$) -- ($(C)+(E)+(F)$);
  \draw ($(C)+(E)+(F)$) -- ($(C)+(D)+(E)+(F)$);
  \draw ($(C)+(E)+(F)$) -- ($(A)+(C)+(E)+(F)$);
  \draw ($(A)+(C)+(F)$) -- ($(A)+(C)+(E)+(F)$);
  \draw ($(A)+(C)+(F)$) -- ($(A)+(B)+(C)+(F)$);
  \draw ($(A)+(B)+(C)$) -- ($(A)+(B)+(C)+(F)$);
  \draw ($(A)+(B)+(C)+(F)$) -- ($(A)+(B)+(C)+(E)+(F)$);
  \draw ($(A)+(C)+(E)+(F)$) -- ($(A)+(B)+(C)+(E)+(F)$);
  \draw ($(A)+(B)+(C)+(E)+(F)$) -- ($(A)+(B)+(C)+(D)+(E)+(F)$);

\end{scope}

\end{tikzpicture}
\caption{\scriptsize\label{fig_doubling}A sequence of doublings, ending with the weak order on $\Sfrak_4$.}
\end{centering}
\end{figure}
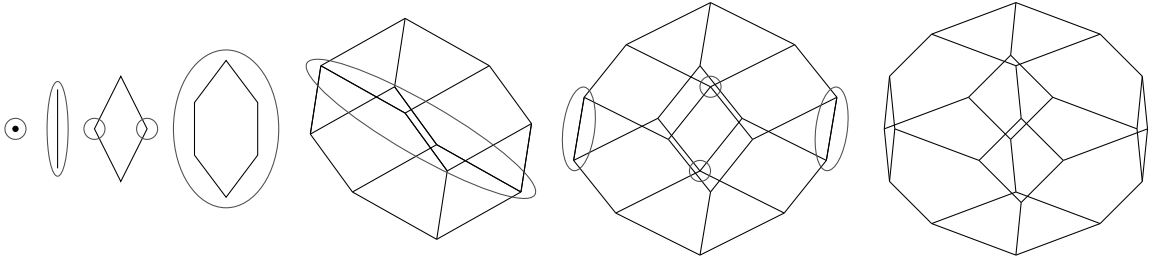

The weak order on permutations is a congruence-uniform lattice; see Figure \ref{fig_doubling} for a sequence of doublings that creates the weak order on $\Sfrak_4$.  The general case is discussed in Example \ref{ex_weak_CN}.

For our purposes, it is easier to employ Reading's characterization of congruence-normal lattices by CN-labelings defined as follows.  For elements $x$ and $y$ of a poset $P$, $y$ \emph{covers} $x$ if $x<y$ and $x\leq z\leq y$ implies $x=z$ or $z=y$ for $z\in P$.  We write $x\lessdot y$ if $y$ covers $x$, and let $\Cov(P)$ denote the set of pairs $(x,y)$ for which $x\lessdot y$.  An \emph{edge-labeling} of a poset $P$ is a function from $\Cov(P)$ to some label set $R$.  Given a lattice $L$ and poset $R$, an edge-labeling $\lambda:\Cov(L)\ra R$ is a \emph{CN-labeling} if $L$ and its dual $L^*$ both satisfy the following condition: For elements $x,y,z\in L$ with $(z,x),(z,y)\in\Cov(L)$ and maximal chains $C_1,C_2\in[z,x\vee y]$ with $x\in C_1,\ y\in C_2$,
\begin{list}{}{}
\item[(CN1)] the elements $x^{\pr}\in C_1,\ y^{\pr}\in C_2$ such that $(x^{\pr},x\vee y),(y^{\pr},x\vee y)\in\Cov(L)$ satisfy
$$\lambda(z,x)=\lambda(y^{\pr},x\vee y),\ \lambda(z,y)=\lambda(x^{\pr},x\vee y);$$
\item[(CN2)] if $(u,v)\in\Cov(C_1)$ with $z<u,\ v<x\vee y$, then $\lambda(z,x)\prec\lambda(u,v)$ and $\lambda(z,y)\prec\lambda(u,v)$; and
\item[(CN3)] the labels on $\Cov(C_1)$ are all distinct.
\end{list}

\begin{theorem}[\cite{reading:lattice}, Theorem 4]\label{thm_CN-labeling}
A finite lattice $L$ is congruence-normal if and only if it admits a CN-labeling.
\end{theorem}

A CN-labeling of the Grassmann-Tamari order $\GT_{3,6}$ is drawn in Figure \ref{fig_P36}.

\begin{theorem}\label{thm_closure_CN}
Let $(S,\prec)$ be a poset with a closure operator.  Assume that
\begin{enumerate}
\item\label{thm_CN_single_step} $\Bic(S)$ is ordered by single-step inclusion,
\item\label{thm_CN_lattice} $W\cup\ov{(X\cup Y)-W}$ is biclosed for $W,X,Y\in\Bic(S)$ with $W\subseteq X\cap Y$, and
\item\label{thm_CN_acyclic} if $x,y,z\in S$ with $z\in\ov{\{x,y\}}-\{x,y\}$ then $x\prec z$ and $y\prec z$.
\end{enumerate}
Then $\Bic(S)$ is a congruence-uniform lattice.
\end{theorem}

\begin{poof}
By Theorem \ref{thm_semidistributive}, we know that $\Bic(S)$ is a semidistributive lattice.  To prove congruence-normality, we verify that $\Bic(S)$ admits a CN-labeling.  Since $\Bic(S)$ is self-dual, the dual conditions will follow from (CN1)-(CN3).

By (\ref{thm_CN_single_step}), we may label a covering relation $X\lessdot Y$ by the unique element in $Y-X$.  These labels are partially ordered by $\prec$.  The property (CN3) is immediate from this definition.

Let $W,X,Y\in\Bic(S)$ such that $X,Y$ both cover $W$.  Let $s,t\in S$ where $X=W\cup\{s\}$ and $Y=W\cup\{t\}$.  By (\ref{thm_CN_lattice}), $X\vee Y=W\cup\ov{\{s,t\}}$ holds, so all of the labels in $[W,X\vee Y]$ lie in $\ov{\{s,t\}}$.  If $C_1$ is a maximal chain in $[X,X\vee Y]$, then the set $X^{\pr}\in C_1$ covered by $X\vee Y$ must be of the form $(X\vee Y)-\{t\}$ as otherwise it would not be biclosed.  Hence (CN1) is satisfied.  Using the relation (\ref{thm_CN_acyclic}), (CN2) is also satisfied.
\end{poof}

\begin{example}\label{ex_weak_CN}
For the closure operator on $\binom{[n]}{2}$ in Example \ref{ex_weak}, we define $\{i,j\}\preceq\{k,l\}$ if $k\leq i<j\leq l$ holds.  By the discussion in Example \ref{ex_weak}, this closure operator satisfies the conditions of Theorem \ref{thm_closure_CN}, so the weak order on permutations is a congruence-uniform lattice.  This holds more generally for the weak order of any finite Coxeter group (\cite[Theorem 6]{caspard.poly-barbut.morvan:cayley} or \cite[Theorem 27]{reading:lattice}).
\end{example}

\section{Biclosed Sets of Segments}\label{sec_segments}

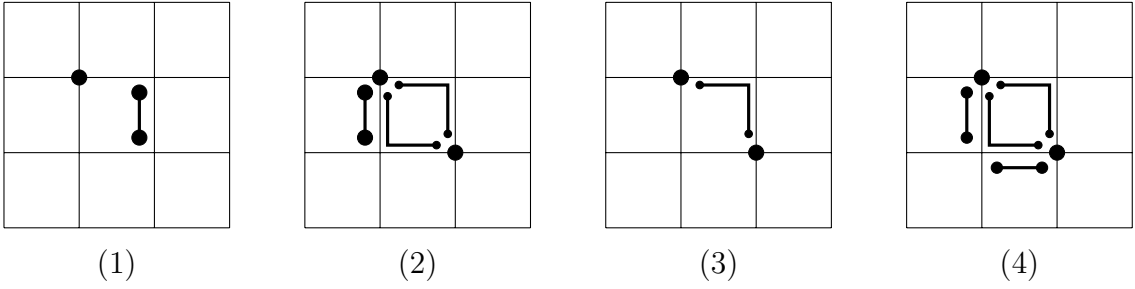
\begin{figure}
\begin{centering}
\begin{tikzpicture}

\begin{scope}

\draw[step=1cm] (0,0) grid (3,3);

\filldraw (1,2) circle(1mm);

\draw[very thick] (1.8,1.8) -- (1.8,1.2);
\filldraw (1.8,1.8) circle(1mm);
\filldraw (1.8,1.2) circle(1mm);

\draw (1.5,-.5) node{(1)};

\end{scope}

\begin{scope}[xshift=4cm]

\draw[step=1cm] (0,0) grid (3,3);

\filldraw (1,2) circle(1mm);
\filldraw (2,1) circle(1mm);
\draw[very thick] (1.25,1.9) -- (1.9,1.9) -- (1.9,1.25);
\filldraw (1.25,1.9) circle(.5mm);
\filldraw (1.9,1.25) circle(.5mm);
\draw[very thick] (1.1,1.75) -- (1.1,1.1) -- (1.75,1.1);
\filldraw (1.1,1.75) circle(.5mm);
\filldraw (1.75,1.1) circle(.5mm);
\draw[very thick] (.8,1.8) -- (.8,1.2);
\filldraw (.8,1.8) circle(1mm);
\filldraw (.8,1.2) circle(1mm);

\draw (1.5,-.5) node{(2)};
\end{scope}

\begin{scope}[xshift=8cm]

\draw[step=1cm] (0,0) grid (3,3);

\filldraw (1,2) circle(1mm);
\filldraw (2,1) circle(1mm);
\draw[very thick] (1.25,1.9) -- (1.9,1.9) -- (1.9,1.25);
\filldraw (1.25,1.9) circle(.5mm);
\filldraw (1.9,1.25) circle(.5mm);

\draw (1.5,-.5) node{(3)};
\end{scope}

\begin{scope}[xshift=12cm]

\draw[step=1cm] (0,0) grid (3,3);

\filldraw (1,2) circle(1mm);
\filldraw (2,1) circle(1mm);
\draw[very thick] (1.25,1.9) -- (1.9,1.9) -- (1.9,1.25);
\filldraw (1.25,1.9) circle(.5mm);
\filldraw (1.9,1.25) circle(.5mm);
\draw[very thick] (1.1,1.75) -- (1.1,1.1) -- (1.75,1.1);
\filldraw (1.1,1.75) circle(.5mm);
\filldraw (1.75,1.1) circle(.5mm);
\draw[very thick] (.8,1.8) -- (.8,1.2);
\filldraw (.8,1.8) circle(.75mm);
\filldraw (.8,1.2) circle(.75mm);
\draw[very thick] (1.2,.8) -- (1.8,.8);
\filldraw (1.2,.8) circle(.75mm);
\filldraw (1.8,.8) circle(.75mm);

\draw (1.5,-.5) node{(4)};
\end{scope}

\end{tikzpicture}
\caption{\scriptsize\label{fig_segments}(1) Two composable segments. (2) A biclosed set $X$ of five segments. (3) $X^{\downarrow}$. (4) $X^{\uparrow}$.}
\end{centering}
\end{figure}

Fix a shape $\lambda$ and let $S$ denote the set of segments supported by $\lambda$.  Two segments $s$ and $t$ are \emph{composable} if $s_{\term}$ is one unit North or West of $t_{\init}$.  If $s$ and $t$ are composable, then the composite $s\circ t$ is the segment containing both $s$ and $t$.  Given a set $X$ of segments of $\lambda$, say $X$ is \emph{closed} if for $s,t\in S$, $s,t\in X$ and $s\circ t\in S$ implies $s\circ t\in X$; see Figure \ref{fig_segments}.  We let $\Bic(S)$ denote the poset of biclosed sets of segments, as in Section \ref{sec_biclosed}.

This closure on segments may be realized as a 2-closure for a certain real vector configuration.  A \emph{cell} of $\lambda$ is a unit square whose four corners are all vertices of $\lambda$.  Let $\Cell(\lambda)$ denote the set of cells of $\lambda$.  To each interval vertex $v$ of $\lambda$, we associate the vector $f_v\in\Rbb^{\Cell(\lambda)}$ where for a cell $c$,
$$f_v(c)=\begin{cases}1 \mbox{ if } v \mbox{ is the SE or NW corner of } c\\-1 \mbox{ if } v \mbox{ is the SW or NE corner of } c\\0 \mbox{ otherwise. }\end{cases}$$

For segments $(v_1,\ldots,v_l)\in S$, set $f_{(v_1,\ldots,v_l)}=\sum_i f_{v_i}$.  It is easy to verify that segments $s,t$ are composable if and only if there exists a segment $u$ such that $f_u=f_s+f_t$.

\begin{example}\label{ex_weak_segments}
Suppose $\lambda$ is a $2\times n$ rectangle.  Labeling the interior vertices $1,\ldots,n-1$ from left to right, a segment $s$ may be identified with the set $\{i,j\}\in\binom{[n]}{2}$ where $i$ is the label on $s_{\init}$ and $j-1$ is the label on $s_{\term}$.  The closure on segments then agrees with the closure on $\binom{[n]}{2}$ defined in Example \ref{ex_weak}.  Hence, $\Bic(S)$ is isomorphic to the weak order on permutations of $[n]$.  Moreover, the vector configuration $\{\frac{1}{\sqrt{2}}f_s:\ s\in S\}$ is the set of positive roots of a root system of type $A_{n-1}$.
\end{example}

\begin{remark}
The vectors $f_s$ for $s\in S$ are called \emph{bending vectors} in \cite{santos.stump.welker:noncrossing}.  Their significance is explained in \cite[Lemma 4.9]{santos.stump.welker:noncrossing}: If $F\stackrel{s}{\ra}F^{\pr}$ are adjacent facets of $\Delta^{NC}_{k,n}$, viewed as a triangulation of the order polytope on a product of chains, then $f_s$ is orthogonal to the ridge $F\cap F^{\pr}$ with $F^{\pr}$ on the positive side.

Given this result, we are led to consider the \emph{bending arrangement} $\Acal_{\lambda}=\{H_s:\ s\in S\}$ where $H_s$ is the hyperplane orthogonal to $f_s$.  This arrangement defines a complete fan on $\Rbb^{k(n-k)}$ whose faces of maximum dimension are called \emph{chambers}.  The \emph{chamber poset} $\Pcal(\Acal_{\lambda})$ is the set of chambers where $c_1\leq c_2$ if $\{s\in S:\ f_s(c_1)>0\}\subseteq\{s\in S:\ f_s(c_2)>0\}$.  Following \cite{reading:lattice_congruence}, we may expect that
\begin{enumerate}
\item\label{arr_lattice} $\Pcal(\Acal_{\lambda})$ is a lattice,
\item\label{arr_quotient} $\GT(\lambda)$ is a lattice quotient of $\Pcal(\Acal_{\lambda})$,
\item\label{arr_fan} $\GT(\lambda)$ is a fan poset on some complete fan $\Fcal$, which is refined by the arrangement fan, and
\item\label{arr_normal} $\Fcal$ is the normal fan of a simple polytope.
\end{enumerate}
However, (\ref{arr_lattice}) is \emph{not} true when $\lambda$ contains a $3\times 3$ square.  The chamber poset naturally injects into the poset of biclosed sets, which we prove is a lattice in Corollary \ref{cor_biclosed_CU}.  Replacing the chamber poset by the poset of biclosed sets, (\ref{arr_quotient}) is a restatement of Corollary \ref{thm_quotient}.  (\ref{arr_fan}) seems to follow from results of \cite{santos.stump.welker:noncrossing}, though we are not sure.  We consider (\ref{arr_normal}) to be an interesting open problem.
\end{remark}

\begin{lemma}\label{lem_SE_corner}
If $X\in\Bic(S)$ and $c$ is a SE-corner of $\lambda$, then $X\setm c$ is a biclosed set of segments of $S\setm c$.
\end{lemma}

\begin{poof}
Let $s,t,u$ be segments supported by $\lambda$ such that $s\circ t=u$.

If $s,t\in X\setm c$ then $u$ is supported by $\lambda\setm c$.  Since $X$ is closed, we conclude $u\in X\setm c$.

If $u\in X\setm c$ then both $s$ and $t$ are supported by $S\setm c$.  Since $S-X$ is closed, either $s$ or $t$ is in $X\setm c$.
\end{poof}

The following description of the closure is immediate from the definition.  We record it here since it is a useful tool in later sections.

\begin{lemma}
For $X\subseteq S$, $\ov{X}$ is the set of segments $s$ such that there exist segments $s_1,\ldots,s_l\in X$ with $s=s_1\circ\cdots\circ s_l$.
\end{lemma}

We partially order $S$ by inclusion; that is, $s\subseteq t$ means $s$ is a subsegment of $t$.

\begin{theorem}\label{thm_bic_seg}
If $\lambda$ is any shape, then
\begin{enumerate}
\item\label{thm_bs_graded} $\Bic(S)$ is ordered by single-step inclusion,
\item\label{thm_bs_closure} $W\cup\ov{(X\cup Y)-W}$ is biclosed for $W,X,Y\in\Bic(S)$ with $W\subseteq X\cap Y$, and
\item\label{thm_bs_acyclic} if $s,t,u\in S$ such that $s\circ t=u$, then $s\subsetneq u$ and $t\subsetneq u$.
\end{enumerate}
\end{theorem}

\begin{poof}
We prove the theorem by induction on the size of $\lambda$.  Let $c$ be a SE-corner of $\lambda$.  We may assume that the vertex $w$ NW of $c$ is an interior vertex as otherwise $S=S\setm c$.

(\ref{thm_bs_graded}): Let $X,Y\in\Bic(S)$ such that $X\subsetneq Y$.  If $s\in Y-X$ is of minimum length, then for any splitting $s=t\circ u$, either $t\in X$ or $u\in X$.  Let $s\in Y-X$ be of maximum length such that for any splitting $s=t\circ u$, either $t\in X$ or $u\in X$.  We prove that $X\cup\{s\}$ is biclosed.

If $X\cup\{s\}$ is not biclosed, then there exists $t\in X$ such that $s\circ t$ or $t\circ s$ is in $S-X$.  Among such segments $t$, choose one of minimum length.  Without loss of generality, we may assume $s\circ t$ is in $S-X$.  Then $s\circ t\in Y$ since $Y$ is closed.  By maximality of $s$, there exists a splitting $s^{\pr}\circ t^{\pr}=s\circ t$ such that $s^{\pr}$ and $t^{\pr}$ are not in $X$.  We distinguish two cases:

(a) Assume $s^{\pr}$ is an initial segment of $s$.  Then $s=s^{\pr}\circ u$ and $t^{\pr}=u\circ t$ for some segment $u$.  By assumption on $s$, we have $u\in X$.  Since $X$ is closed, this forces $t^{\pr}\in X$, a contradiction.

(b) Assume $s$ is an initial segment of $s^{\pr}$.  Then $s^{\pr}=s\circ u$ and $t=u\circ t^{\pr}$ for some segment $u$.  Since $t\in X,\ t^{\pr}\notin X$, we deduce $u\in X$.  Since $u$ is shorter than $t$, we deduce $s^{\pr}\in X$, a contradiction.

Hence, $X\cup\{s\}$ is biclosed.

(\ref{thm_bs_closure}): Let $W\in\Bic(S)$.  Assume, for $W^{\pr}\in\Bic(S)$ with $W\subsetneq W^{\pr}$:

$$W^{\pr}\cup\ov{(X\cup Y)-W^{\pr}}\mbox{ is biclosed for }X,Y\in\Bic(S)\mbox{ with }W^{\pr}\subseteq X\cap Y.$$

Let $X,Y\in\Bic(S)$ such that $W\subseteq X\cap Y$.  We may assume that $W$ is a maximal biclosed set contained in $X\cap Y$.  If $X\subseteq Y$, the result is immediate.  If $X$ and $Y$ are incomparable, then by (\ref{thm_bs_graded}), there exists $s\in X-W,\ t\in Y-W$ such that $W\cup\{s\}$ and $W\cup\{t\}$ are biclosed.

Set $Z=W\cup\ov{\{s,t\}}$.  If $s$ and $t$ are not composable, then $Z=W\cup\{s,t\}$ is biclosed.  If $s\circ t=u$, we claim that $Z=W\cup\{s,t,u\}$ is biclosed.

If $Z$ is not closed, then there exists $v\in W$ such that $v\circ u$ or $u\circ v$ is in $S-W$.  We may assume without loss of generality that $v\circ u$ is in $S-W$.  Since $W\cup\{s\}$ is closed and $v\circ s\in S$, we have $v\circ s\in W$.  But $W\cup\{t\}$ is closed, so $v\circ s\circ t$ is in $W$, a contradiction.  Hence, $Z$ is closed.

If $S-Z$ is not closed, then there exists a splitting $s^{\pr}\circ t^{\pr}=u$ such that $s^{\pr}$ and $t^{\pr}$ are in $S-Z$.  Then either $s$ is an initial subsegment of $s^{\pr}$ or $t$ is a terminal subsegment of $t^{\pr}$.  Without loss of generality, we may assume $s$ is an initial subsegment of $s^{\pr}$.  Then there exists a segment $u^{\pr}$ with $s\circ u^{\pr}=s^{\pr}$ and $u^{\pr}\circ t^{\pr}=t$.  Since $W\cup\{s\}$ is closed, the condition $s\circ u^{\pr}=s^{\pr}$ implies $u^{\pr}\notin W$.  However, as $S-(W\cup\{t\})$ is closed, the latter condition implies $u^{\pr}\in W$, a contradiction.

Therefore, $Z$ is biclosed.  Applying the assumption with $W^{\pr}=W\cup\{s\}$, we deduce that
$$W\cup\{s\}\cup\ov{(X\cup Z)-(W\cup\{s\})}$$
is biclosed.  Similarly,
$$W\cup\{t\}\cup\ov{(Y\cup Z)-(W\cup\{t\})}$$
is biclosed.  As both of these sets contain $Z$, we deduce that
$$Z\cup\ov{((W\cup\{s\}\cup\ov{(X\cup Z)-(W\cup\{s\})})\cup(W\cup\{t\}\cup\ov{(Y\cup Z)-(W\cup\{t\})}))- Z}$$
is biclosed.  This set is equal to
$$Z\cup\ov{\ov{(X\cup Y\cup Z)- W}- Z}.$$
But,
$$X\cup Y\subseteq Z\cup\ov{\ov{(X\cup Y\cup Z)- W}- Z}\subseteq W\cup\ov{(X\cup Y\cup Z)- W}=W\cup\ov{(X\cup Y)- W}\subseteq\ov{X\cup Y}.$$
Since $\ov{X\cup Y}$ is the smallest closed set containing $X\cup Y$, we deduce the equality
$$Z\cup\ov{\ov{(X\cup Y\cup Z)- W}- Z}=\ov{X\cup Y}.$$
Hence, $W\cup\ov{(X\cup Y)- W}$ is biclosed, as desired.

(\ref{thm_bs_acyclic}): This is immediate from the definitions.
\end{poof}

Applying Theorem \ref{thm_closure_CN}, we deduce

\begin{corollary}\label{cor_biclosed_CU}
$\Bic(S)$ is a congruence-uniform lattice.
\end{corollary}

\begin{remark}
The hypotheses of Theorems \ref{thm_semidistributive} and \ref{thm_closure_CN} were chosen with two examples in mind, namely the 2-closure on finite root systems and the closure operator defined in this section.  For the 2-closure on a real simplicial hyperplane arrangement, the first two hypotheses hold, but the third may not.  In this case, a weaker version of the acyclic condition is enough to prove congruence-normality \cite[Theorem 25]{reading:lattice}.
\end{remark}

\section{A quotient of $\Bic(S)$}\label{sec_quotient}

Given a biclosed set $X$ of segments, let $X^{\downarrow}$ be the set of segments $s$ in $X$ such that $t$ is in $X$ whenever $t$ is a SW-subsegment of $s$.  Let $X^{\uparrow}$ be the set of segments $s$ such that there exists $t$ in $X$ that is a NE-subsegment of $s$.  An example is shown in Figure \ref{fig_segments}.

Transposition of shapes $\lambda\ra\lambda^{\tr}$ induces a map on segments $s\mapsto s^{\tr}$.  Given a set $X$ of segments of $\lambda$, we let $X^{\tr}$ denote the set of transposed segments of $\lambda^{\tr}$.  Transposition commutes with complementation.  Let $X^{c\tr}$ be the composition of these two involutions.

\begin{claim}\label{claim_complement}
For $X\subseteq S$, 
$$(X^{\uparrow})^{c\tr}=(X^{c\tr})^{\downarrow}.$$
\end{claim}

\begin{poof}
A segment $s$ is in $(X^{\uparrow})^{c\tr}$ if and only if $s^{\tr}$ is not a segment in $X^{\uparrow}$.  But this holds exactly when none of the NE-subsegments of $s^{\tr}$ are in $X$.  This occurs if none of the SW-subsegments of $s$ are in $X^{\tr}$, which is equivalent to $s\in(X^{c\tr})^{\downarrow}$.
\end{poof}

\begin{claim}
If $X$ is biclosed, then
\begin{enumerate}
\item $X^{\tr}$ is biclosed,
\item $X^c$ is biclosed,
\item $X^{\downarrow}$ is biclosed, and
\item $X^{\uparrow}$ is biclosed.
\end{enumerate}
\end{claim}

\begin{poof}
Parts (1) and (2) are immediate from the definitions.  Part (4) follows from (1)-(3) with Claim \ref{claim_complement}.  We verify part (3).

Let $s,t\in X^{\downarrow}$ such that $s\circ t$ is a segment.  Since $X$ is biclosed, $s\circ t$ is in $X$.  If $u$ is a SW-subsegment of $s\circ t$, then either $u\subseteq s,\ u\subseteq t$ or neither inequality holds.  In the first two cases, it follows that $u\in X$ from $s,t\in X^{\downarrow}$.  In the remaining case, we may divide $u$ into two pieces $u=u_1\circ u_2$ where $u_1$ is a SW-subsegment of $s$ and $u_2$ is a SW-subsegment of $t$.  Hence, $u_1,u_2\in X$, so also $u\in X$.  Therefore, $X^{\downarrow}$ is closed.

On the other hand, if $s\in X^{\downarrow}$ such that $s=t\circ u$, then either $t$ or $u$ is a SW-subsegment.  Hence, $X^{\downarrow}$ is co-closed as well.
\end{poof}

\begin{claim}\label{claim_idem_order}
The maps $X\mapsto X^{\downarrow}$ and $X\mapsto X^{\uparrow}$ are idempotent and order-preserving.
\end{claim}

\begin{poof}
The order-preserving assertion is immediate from the definition.  It remains to prove the maps are idempotent.

For segments $s,t,u$, if $s$ is a SW-subsegment of $t$ and $t$ is a SW-subsegment of $u$, then $s$ is a SW-subsegment of $u$.  Hence, for $u\in X$, if every subsegment of $u$ is in $X$, then every subsegment of $u$ is also in $X^{\downarrow}$.  The claim follows immediately.
\end{poof}

\begin{claim}
$(X^{\downarrow})^{\uparrow}=X^{\uparrow}$.  Dually, $(X^{\uparrow})^{\downarrow}=X^{\downarrow}$.
\end{claim}

\begin{poof}
The forward inclusion $(X^{\downarrow})^{\uparrow}\subseteq X^{\uparrow}$ follows from Claim \ref{claim_idem_order}.

If $s\in X^{\uparrow}$, then there exists $t_0\in X$ such that $t_0$ is a NE-subsegment of $s$.  If $t_0\notin X^{\downarrow}$, then there exists a SW-subsegment $u_0$ of $t_0$ that is not in $X$.  Then $t_0=u_0^{\pr}\circ u\circ u_0^{\prpr}$ where $u_0^{\pr}$ and $u_0^{\prpr}$ are (possibly empty) NE-subsegments of $t_0$.  Since $X$ is biclosed, either $u_0^{\pr}$ or $u_0^{\prpr}$ is in $X$.  In particular, $t_0$ has a NE-subsegment $t_1$ that is in $X$.  Continuing in this manner, we produce a segment $t\in X^{\downarrow}$ that is a NE-subsegment of $s$.  Hence, $s\in(X^{\downarrow})^{\uparrow}$, as desired.

The second claim follows from the first via Claim \ref{claim_complement}.
\end{poof}

Let $\Theta$ be the equivalence relation on $\Bic(S)$ where $X\equiv Y\mod\Theta$ if $X^{\downarrow}=Y^{\downarrow}$.

\begin{theorem}\label{thm_congruence}
$\Theta$ is a lattice congruence on $\Bic(S)$.
\end{theorem}

\begin{poof}
We prove that $\Theta$ is a lattice congruence following Proposition \ref{prop_congruence}.

Let $X,Y\in\Bic(S)$ such that $X\equiv Y\mod\Theta$.  Since $X^{\downarrow}\subseteq X$, it follows that $X^{\downarrow}$ is the smallest element in $[X]$.  Since
$$X^{\uparrow}=(X^{\downarrow})^{\uparrow}=(Y^{\downarrow})^{\uparrow}=Y^{\uparrow},$$
it follows that $X^{\uparrow}$ is the largest element in $[X]$.  If $Z\in[X^{\downarrow},X^{\uparrow}]$, then
$$X^{\downarrow}=(X^{\downarrow})^{\downarrow}\subseteq Z^{\downarrow}\subseteq (X^{\uparrow})^{\downarrow}=X^{\downarrow},$$
so the interval $[X^{\downarrow},X^{\uparrow}]$ is the equivalence class of $X$.

The maps $\pi_{\downarrow}(X)=X^{\downarrow}$ and $\pi^{\uparrow}(X)=X^{\uparrow}$ are order-preserving by Claim \ref{claim_idem_order}.
\end{poof}

\begin{example}\label{ex_312}
Let $\lambda$ be the $2\times n$ rectangle from Example \ref{ex_weak_segments}.  If $X$ is a biclosed subset of $S$, then $X^{\downarrow}$ is the set obtained by removing horizontal segments for which some initial part is not in $X$.  The set $X^{\uparrow}$ is obtained by adding horizontal segments to $X$ for which some initial part is not in $X$ but the corresponding terminal part is in $X$.  By this observation it follows that $X^{\uparrow}$ is the largest biclosed set for which $(X^{\uparrow})^{\downarrow}=X^{\downarrow}$.  In particular, the equivalence classes are all closed intervals of the form $[X^{\downarrow},X^{\uparrow}]$ for some $X\in\Bic(S)$.  Moreover, $\pi^{\uparrow}(X)=X^{\uparrow}$ and $\pi_{\downarrow}(X)=X^{\downarrow}$, so $\pi^{\uparrow}$ and $\pi_{\downarrow}$ are both order-preserving maps, thus verifying Theorem \ref{thm_congruence} in this case.  The argument for general shapes follows similar reasoning.

When $\lambda$ is a $2\times n$ rectangle, the bijection in Example \ref{ex_weak_segments} takes biclosed sets $X$ for which $X^{\downarrow}=X$ to inversion sets of 312-avoiding permutations.  Indeed, if a permutation $\sigma=\sigma_1\cdots\sigma_n$ contains a 312 pattern, say with values $i<j<k$, then the corresponding biclosed set $X$ has a long segment labeled $\{i,k\}$ for which the initial part $\{i,j\}$ is not in $X$.
\end{example}

\section{Proof of Theorem \ref{thm_main}}\label{sec_tamari}

For this section, we fix a shape $\lambda$, and let $S$ denote the set of segments of $\lambda$.  Furthermore, we let $E_V$ denote the set of interior vertical edges in $\lambda$ and let $\Pcal$ be the set of paths supported by $\lambda$.

\begin{figure}
\begin{centering}
\begin{tikzpicture}

\begin{scope}

\draw[step=1cm] (0,0) grid (3,3);

\filldraw (1,2) circle(1mm);
\filldraw (2,1) circle(1mm);
\filldraw (1,1) circle(1mm);

\draw[very thick] (1.1,1.75) -- (1.1,1.1) -- (1.75,1.1);
\filldraw (1.1,1.75) circle(.5mm);
\filldraw (1.75,1.1) circle(.5mm);

\draw[very thick] (.8,1.8) -- (.8,1.2);
\filldraw (.8,1.8) circle(1mm);
\filldraw (.8,1.2) circle(1mm);

\draw[very thick] (1.2,.8) -- (1.8,.8);
\filldraw (1.2,.8) circle(1mm);
\filldraw (1.8,.8) circle(1mm);

\draw[decorate,decoration=zigzag,color=orange] (1,3) -- (1,2) -- (2,2) -- (2,1) -- (2,0);

\draw[purple, very thick] (.8,2.5) -- (1.2,2.5);

\end{scope}

\begin{scope}[xshift=4cm]

\draw[step=1cm] (0,0) grid (3,3);

\filldraw (1,2) circle(1mm);
\filldraw (2,1) circle(1mm);
\filldraw (1,1) circle(1mm);

\draw[very thick] (1.1,1.75) -- (1.1,1.1) -- (1.75,1.1);
\filldraw (1.1,1.75) circle(.5mm);
\filldraw (1.75,1.1) circle(.5mm);

\draw[very thick] (.8,1.8) -- (.8,1.2);
\filldraw (.8,1.8) circle(1mm);
\filldraw (.8,1.2) circle(1mm);

\draw[very thick] (1.2,.8) -- (1.8,.8);
\filldraw (1.2,.8) circle(1mm);
\filldraw (1.8,.8) circle(1mm);

\draw[decorate,decoration={coil,segment length=4pt},color=green] (1,3) -- (1,1) -- (3,1);

\draw[purple, very thick] (.8,1.5) -- (1.2,1.5);

\end{scope}

\begin{scope}[xshift=8cm]

\draw[step=1cm] (0,0) grid (3,3);

\filldraw (1,2) circle(1mm);
\filldraw (2,1) circle(1mm);
\filldraw (1,1) circle(1mm);

\draw[very thick] (1.1,1.75) -- (1.1,1.1) -- (1.75,1.1);
\filldraw (1.1,1.75) circle(.5mm);
\filldraw (1.75,1.1) circle(.5mm);

\draw[very thick] (.8,1.8) -- (.8,1.2);
\filldraw (.8,1.8) circle(1mm);
\filldraw (.8,1.2) circle(1mm);

\draw[very thick] (1.2,.8) -- (1.8,.8);
\filldraw (1.2,.8) circle(1mm);
\filldraw (1.8,.8) circle(1mm);

\draw[decorate,decoration=zigzag,color=blue] (0,2) -- (1,2) -- (2,2) -- (2,1) -- (3,1);

\draw[purple, very thick] (1.8,1.5) -- (2.2,1.5);

\end{scope}

\begin{scope}[xshift=12cm]

\draw[step=1cm] (0,0) grid (3,3);

\filldraw (1,2) circle(1mm);
\filldraw (2,1) circle(1mm);
\filldraw (1,1) circle(1mm);

\draw[very thick] (1.1,1.75) -- (1.1,1.1) -- (1.75,1.1);
\filldraw (1.1,1.75) circle(.5mm);
\filldraw (1.75,1.1) circle(.5mm);

\draw[very thick] (.8,1.8) -- (.8,1.2);
\filldraw (.8,1.8) circle(1mm);
\filldraw (.8,1.2) circle(1mm);

\draw[very thick] (1.2,.8) -- (1.8,.8);
\filldraw (1.2,.8) circle(1mm);
\filldraw (1.8,.8) circle(1mm);

\draw[decorate,decoration={coil,segment length=4pt},color=red] (2,0) -- (2,2) -- (0,2);

\draw[purple, very thick] (1.8,.5) -- (2.2,.5);

\end{scope}

\end{tikzpicture}
\caption{\scriptsize\label{fig_eta}The four non-vertical and non-horizontal paths in $\eta(X)$ where $X$ is the set of black segments.}
\end{centering}
\end{figure}

We define a function $\eta:\Bic(S)\ra 2^{\Pcal}$ as follows.  Let $X\in\Bic(S)$ be given.  If $e\in E_V$ is an edge from $u$ to $v$, let $p_e$ be the path such that for interior vertices $u^{\pr}\in p_e[\cdot,u]$ and $v^{\pr}\in p_e[v,\cdot]$:
\begin{list}{}{}
\item[(i)] if $p_e[u^{\pr},u]$ is (not) in $X$ then $p_e$ enters $u^{\pr}$ from the North (West); and
\item[(ii)] if $p_e[v,v^{\pr}]$ is (not) in $X$ then $p_e$ leaves $v^{\pr}$ to the East (South).
\end{list}

Let $\eta(X)$ be the union of $\{p_e:\ e\in E_V\}$ with the set of horizontal paths supported by $\lambda$.

\begin{example}
If $X$ is the biclosed set of six black segments in Figure \ref{fig_eta}, each of the six interior vertical edges corresponds to a non-horizontal path in $\eta(X)$.  In Figure \ref{fig_eta}, the four paths corresponding to the four marked purple edges are drawn.  The other two vertical edges correspond to vertical paths.  This is the same collection of paths as in Example \ref{ex_facet}.
\end{example}

\begin{claim}\label{claim_eta_kiss}
Let $p$ be a path in $\eta(X)$ containing a segment $s$.
\begin{enumerate}
\item\label{claim_ek_1} If $p$ enters $s$ from the West and leaves $s$ to the South, then $s$ is not in $X$.
\item\label{claim_ek_2} Similarly, if $p$ enters $s$ from the North and leaves $s$ to the East, then $s$ is in $X$.
\end{enumerate}
\end{claim}

\begin{poof}
We prove \ref{claim_ek_1}.  The proof of \ref{claim_ek_2} is similar.

Let $e$ be the interior vertical edge with $p=p_e$.  Assume $p$ enters $s$ from the West and leaves $s$ to the South.  We prove that $s$ is not in $X$ by considering several cases.

If $e$ is contained in $s$, then by construction $s[\cdot,e_{\init}]$ and $s[e_{\term},\cdot]$ are not in $X$.  If $e$ precedes $s$ in $p$, then $p[e_{\term},s_{\term}]$ is not in $X$ while $p[e_{\term},s_{\term}]-s$ is in $X$.  Finally, if $e$ comes after $s$ in $p$, then $p[s_{\init},e_{\init}]$ is not in $X$ while $p[s_{\init},e_{\init}]-s$ is in $X$.  In each case, we conclude that $s$ is not in $X$ since $X$ is biclosed.
\end{poof}

\begin{claim}\label{claim_eta_max}
If $e$ and $e^{\pr}$ are distinct interior vertical edges, then $p_e$ and $p_{e^{\pr}}$ are distinct paths.
\end{claim}

\begin{poof}
Assume to the contrary that $p_e$ and $p_{e^{\pr}}$ are the same.  Without loss of generality, we may assume that $e$ precedes $e^{\pr}$ in $p_e$.  By definition of $p_e$, the segment $p_e[e_{\term},e^{\pr}_{\init}]$ is not in $X$.  By the definition of $p_{e^{\pr}}$, the same segment $p_{e^{\pr}}[e_{\term},e^{\pr}_{\init}]$ is in $X$, a contradiction.
\end{poof}

\begin{proposition}\label{prop_eta_def}
$\eta(X)$ is a maximal collection of non-kissing paths.
\end{proposition}

\begin{poof}
Suppose that $\eta(X)$ contains two paths $p_{e_1},p_{e_2}$ kissing along a common segment $s$.  By Claim \ref{claim_eta_kiss}, $s$ must be both in $X$ and not in $X$, a contradiction.  Hence, $\eta(X)$ is a set of non-kissing paths.

By Claim \ref{claim_eta_max}, $\eta(X)$ is of maximal size.
\end{poof}

By Proposition \ref{prop_eta_def}, $\eta$ is a map from $\Bic(S)$ to $\Fcal(\Delta^{NK}(\lambda))$.

\begin{claim}\label{claim_eta_top}
$p_e$ is the top path at $e$ in $\eta(X)$ (i.e. $p_e$ is maximum with respect to the total order $\prec_e$ from Section \ref{sec_non-kissing}).
\end{claim}

\begin{poof}
Let $e^{\pr}$ be an edge distinct from $e$ such that $p_{e^{\pr}}$ contains $e$.  Without loss of generality, we may assume $e$ precedes $e^{\pr}$ in $p_{e^{\pr}}$.  By definition of $p_{e^{\pr}}$, $X$ contains the segment $p_{e^{\pr}}[e_{\term},e^{\pr}_{\init}]$.

Let $s$ be the initial segment of $p_e[e_{\term},\cdot]$ along which $p_e$ and $p_{e^{\pr}}$ agree.  Assume $p_e$ leaves $s$ to the South and $p_{e^{\pr}}$ leaves to the East.  By definition of $p_e$, this implies $s$ is not in $X$.  By definition of $p_{e^{\pr}}$, the segment $p_{e^{\pr}}[e_{\term},e^{\pr}_{\init}]-s$ is also not in $X$.  As $X$ is co-closed, this implies $p_{e^{\pr}}[e_{\term},e^{\pr}_{\init}]$ is not in $X$, a contradiction.  Therefore, $p_e$ is the top path at $e$.
\end{poof}

Let $F$ be a facet of $\Delta^{NK}(\lambda)$.  For each path $p\in F$, let $A_p$ be the set of SW-subsegments of $p$.  That is, $A_p$ consists of the segments $p[v,v^{\pr}]$ such that $p$ enters $v$ from the North and leaves $v^{\pr}$ to the East.  Set $\ds \phi(F)=\ov{\bigcup_{p\in F}A_p}$.  A priori, $\phi$ is a map from facets of $\Delta^{NK}(S)$ to sets of segments.

\begin{claim}\label{claim_phi_def}
$\phi(F)$ is a biclosed set of segments.
\end{claim}

\begin{poof}
It is sufficient to show that $A_p$ is a biclosed set as this would imply
$$\bigvee_{p\in F} A_p=\ov{\bigcup_{p\in F}A_p}.$$

No two segments in $A_p$ are composable, so it is closed.  Let $s\in A_p$, and let $t,u$ be segments such that $t\circ u=s$.  If the edge separating $t$ and $u$ is horizontal, then $t$ is in $A_p$.  If the edge separating $t$ and $u$ is vertical, then $u$ is in $A_p$.  Hence, $A_p$ is biclosed.
\end{poof}

Now we have defined functions $\eta:\Bic(S)\ra\Fcal(\Delta^{NK}(\lambda))$ and $\phi:\Fcal(\Delta^{NK}(\lambda))\ra\Bic(S)$.  We next show that $\eta$ is surjective.

\begin{claim}\label{claim_eta_surjective}
The composite $\eta\circ\phi$ is equal to the identity on $\Fcal(\Delta^{NK}(\lambda))$.
\end{claim}

\begin{poof}
Let $F\in\Fcal(\Delta^{NK}(\lambda))$.  Given a non-horizontal path $p$ in $F$, we prove that $p$ is in $\eta(\phi(F))$.  Suppose $p$ is on top at edge $e$.  Let $q$ be the path associated to $e$ in $\eta(\phi(F))$.  If $q=p$, we are done.  Otherwise, we may assume that $p$ and $q$ are distinct after $e$.  If not, then a similar argument may be used when $p$ and $q$ are distinct before $e$.

Let $s$ be the longest segment along which $p$ and $q$ agree starting from $e_{\term}$.  If $p$ leaves $s$ to the East, then $s$ is in $A_p$.  This forces $q$ to leave $s$ to the East as well, contradicting the maximality of $s$.  Hence, $q$ leaves $s$ to the East and $p$ leaves $s$ to the South.

Since $s$ is in $\ov{\bigcup_{p\in F}A_p}$, the segment may be decomposed as $s=s_1\circ\cdots\circ s_l$ where $s_i\in A_{p_i}$ for $i\in[l]$.  If $s_1\notin A_p$, then $p<_e p_1$, a contradiction.  Let $k$ be the smallest index for which $s_1\circ\cdots\circ s_k\notin A_p$.  Then $p$ enters $s_k$ from the West and leaves to the South.  But $s_k\in A_{p_k}$, so $p$ and $p_k$ are kissing, a contradiction.

But this means $s\in A_p$, in contradiction to the assumption that $p$ leaves $s$ to the South.  Hence, we conclude $p=q$.
\end{poof}

\begin{claim}\label{claim_phi_down}
For $X\in\Bic(S)$, $\phi\circ\eta(X)=X^{\downarrow}$.
\end{claim}

\begin{poof}
We prove
\begin{list}{}{}
\item[(a)] $X^{\downarrow}\subseteq\phi\circ\eta(X)$
\item[(b)] $\phi\circ\eta(X)\subseteq X$
\item[(c)] For $s\in\phi\circ\eta(X)$, if $s^{\pr}$ is a SW-subsegment of $s$, then $s^{\pr}\in\phi\circ\eta(X)$.
\end{list}

The claim is then immediate.

(a) Let $s\in X^{\downarrow}$.  We show $s$ is in $\phi\circ\eta(X)$ by induction on the length of $s$.  Let $e$ be the vertical edge with terminal vertex $s_{\init}$.

Let $t$ be the initial subsegment of $s$ that coincides with $p_e$ in $\eta(X)$.  If $s=t$, then $p_e$ leaves $s$ to the East, so $s$ is in $A_{p_e}$.  Assume $s=t\circ t^{\pr}$ for some segment $t^{\pr}$.  If $s$ leaves $t$ to the East, then $t$ is in $X$, so $p_e$ also leaves to the East.  Hence, $s$ leaves $t$ to the South while $p_e$ leaves to the East.  By the induction hypothesis, $t^{\pr}$ is in $\phi\circ\eta(X)$.  Hence, $s$ is in $\phi\circ\eta(X)$, as desired.

(b) Let $p\in\eta(X)$.  It suffices to show that $A_p$ is a subset of $X$.  Suppose $p$ is on top at some vertical edge $e$.  By Claim \ref{claim_eta_top}, $p=p_e$ in the construction of $\eta(X)$.  Fix $s\in A_p$.  If $s$ contains $e$, then since $s$ is a SW-subsegment of $p$, $X$ contains both $s[\cdot,e_{\init}]$ and $s[e_{\term},\cdot]$.  As $X$ is closed, this implies $s\in X$.  If $e$ precedes $s$, then $X$ contains $p[e_{\term},s_{\term}]$ but not $p[e_{\term},s_{\term}]-s$.  As $X$ is co-closed, this implies $s\in X$.  If $e$ appears after $s$, then $X$ contains $p[s_{\init},e_{\init}]$ but not $p[s_{\init},e_{\init}]-s$.  As $X$ is co-closed, this forces $s\in X$.  Therefore, we conclude that $A_p\subseteq X$ holds.

(c) Let $s\in\phi\circ\eta(X)$ and let $s^{\pr}$ be a SW-subsegment of $s$.  We show that $s^{\pr}$ is in $X$.  Then $s=s_1\circ\cdots\circ s_l$ for some segments $s_i\in A_{p_i}$ and paths $p_i\in\eta(X)$.  As $s^{\pr}$ is a subsegment of $s$, there exist indices $i\leq j$ for which $s^{\pr}=s_i^{\pr}\circ s_{i+1}\circ\cdots\circ s_{j-1}\circ s_j^{\pr}$ where $s_i^{\pr}$ is a subsegment of $s_i$ and $s_j^{\pr}$ is a subsegment of $s_j$.  Since $s^{\pr}$ is a SW-subsegment of $s$, it follows that $s_i^{\pr}$ is a SW-subsegment of $s_i$ and $s_j^{\pr}$ is a SW-subsegment of $s_j$.  Hence, $s_i^{\pr}\in A_{p_i}$ and $s_j^{\pr}\in A_{p_j}$.  We conclude that $s^{\pr}$ is in $\phi\circ\eta(X)$.
\end{poof}

Using Claim \ref{claim_phi_down}, it is easy to show that the fibers of $\eta$ are equivalence classes of $\Theta$.

\begin{claim}\label{claim_quotient}
For $X,Y\in\Bic(S)$, $\eta(X)=\eta(Y)$ if and only if $X^{\downarrow}=Y^{\downarrow}$.
\end{claim}

\begin{poof}
Assume $\eta(X)=\eta(Y)$.  Then
$$X^{\downarrow}=\phi\circ\eta(X)=\phi\circ\eta(Y)=Y^{\downarrow}.$$

Now assume $X^{\downarrow}=Y^{\downarrow}$.  Then
$$\eta(X)=\eta\circ\phi\circ\eta(X)=\eta(X^{\downarrow})=\eta(Y^{\downarrow})=\eta\circ\phi\circ\eta(Y)=\eta(Y),$$
as desired.
\end{poof}

\begin{claim}\label{claim_covers}
For $F\in\Fcal(\Delta^{NK}(\lambda))$,
$$\{s\in S:\ \exists F^{\pr}\in\Fcal(\Delta^{NK}(\lambda)),\ F^{\pr}\stackrel{s}{\ra}F\}=\{s\in\phi(F):\ \phi(F)-\{s\}\ \mbox{is biclosed}\}.$$
Moreover, for adjacent facets $F,F^{\pr}$, if $F^{\pr}\stackrel{s}{\ra}F$, then $\eta(\phi(F)-\{s\})=F^{\pr}$.
\end{claim}

\begin{poof}
We first show the forward inclusion.  Let $s\in S$ and $F^{\pr}\in\Fcal(\Delta^{NK}(\lambda))$ adjacent to $F$ with $F^{\pr}\stackrel{s}{\ra}F$.  Let $p\in F-F^{\pr},\ p^{\pr}\in F^{\pr}-F$, so $p$ and $p^{\pr}$ kiss along $s$.

Assume $\phi(F)-\{s\}$ is not biclosed.  If $\phi(F)-\{s\}$ is not closed, then $s=s_1\circ\cdots\circ s_l$ for some $s_i\in A_{p_i}$ and $p_i\in F$ with $l>1$.  If $s_1\notin A_p$, then $p_1$ and $p^{\pr}$ are kissing, an impossibility.  Let $k\geq 1$ and suppose $s_1\circ\cdots\circ s_k$ is in $A_p$ but not $s_1\circ\cdots\circ s_{k+1}$.  Then $p_{k+1}$ and $p^{\pr}$ are kissing, which is again impossible.  Hence, $s_1\circ\cdots\circ s_{l-1}$ is in $A_p$.  Since $p^{\pr}$ leaves $s_l$ to the South and enters from the West, $p^{\pr}$ and $s_l$ are kissing, a contradiction.  Hence, $\phi(F)-\{s\}$ is closed.

Now assume $\phi(F)-\{s\}$ is not co-closed.  We may assume that there exists $t\in S-\phi(F)$ such that $s\circ t\in\phi(F)$ but $t\notin\phi(F)$.  We prove $t\in\phi(F)$ by induction on the length of $t$, which gives a contradiction.

Then $s\circ t=u_1\circ\cdots\circ u_l$ where $u_i\in A_{p_i}$ and $p_i\in F$ for all $i$. If $p$ leaves $u_1$ to the South, then $p^{\pr}$ and $p_1$ are kissing, an impossibility.  As before, we determine that $s\circ t$ is a SW-subsegment of $p$.  As $s$ is a SW-subsegment of $p$, this implies that the edge $e$ between $s$ and $t$ is horizontal.  Moreover, $p$ is the bottom path at $e$ in $F$.  Let $e_1$ be the horizontal edge after $t_{\term}$, and let $q_1$ be the bottom path at $e_1$.  Then $q_1$ and $p$ agree along a terminal subsegment $t_1$ of $t$, where $p$ enters $t_1$ from the West and $q_1$ enters $t_1$ from the North.  Moreover, $s\circ(t-t_1)\in\phi(F)$, so $t-t_1\in\phi(F)$ by induction.  But $t_1\in A_{q_1}$, so $t=(t-t_1)\circ t_1\in\phi(F)$.

Now we prove the reverse inclusion.  Let $s\in\phi(F)$ such that $\phi(F)-\{s\}$ is biclosed.  We prove that $\eta(\phi(F)-\{s\})$ is adjacent to $F$ and $\eta(\phi(F)-\{s\})\stackrel{s}{\ra}F$.

Since $\phi(F)=\phi(\eta\circ\phi(F))$ and $\phi(F)-\{s\}$ is co-closed, there exists a path $p$ such that $s\in A_p$.  Let $e$ be the vertical edge above $s_{\init}$.  Since $s\in A_p$, it follows that $p$ contains $e$.

We show that $p$ is the top path at $e$.  Suppose not, and let $q$ be the top path at $e$.  If $q$ does not contain $s$, then let $t$ be the largest subsegment of $s$ along which $p$ and $q$ agree.  Then $q$ leaves $t$ to the East and $p$ leaves $t$ to the South, so $t\in A_{q}$ and $s-t\in A_p$.  But this is impossible since $\phi(F)-\{s\}$ is co-closed.

Suppose $q$ contains $s$ and let $v$ be the first vertex after $s$ such that $q$ leaves $v$ to the East and $p$ leaves $v$ to the South.  (We note that if $q$ and $p$ agree after $s$, then we may deduce a contradiction in a similar way where we take $v$ to occur before $s$.)  Let $t$ be the segment $p[s_{\init},v]$.  Since $t\in A_{q}$ and $\phi(F)-\{s\}$ is co-closed, $t-s$ is in $\phi(F)$.  Let $u_1,\ldots,u_l$ be segments such that $u_i\in A_{p_i}$ for some paths $p_i$ and $t-s=u_1\circ\cdots\circ u_l$.  Since $p$ and $p_1$ are non-kissing, $p$ must leave $u_1$ to the East.  Similarly, we deduce that $p$ leaves $u_2,\ldots,u_l$ to the East.  But $p$ leaves $u_l$ to the South, a contradiction.

We have now determined that $p=q$.  As $p$ was chosen as an arbitrary path containing $s$ as a SW-subsegment, it follows that $p$ is the unique such path.

Let $e$ be the edge above $s_{\init}$ and $e^{\pr}$ the edge below $s_{\term}$.  By definition, $\eta(\phi(F)-\{s\})$ differs from $F$ by at most two paths, namely $p_e$ and $p_{e^{\pr}}$.  We claim $p_e$ is in $F$ and it is the top path in $F$ at $e^{\pr}$.  As $p\notin\eta(\phi(F)-\{s\})$, it would follow that $F$ and $\eta(\phi(F)-\{s\})$ are adjacent facets and the new path in $\eta(\phi(F)-\{s\})$ kisses $p$ along $s$.

Let $q$ be the top path in $F$ at $e^{\pr}$.  Since $\phi(F)$ contains all SW-subsegments of $s$ and $p$ is the top path in $F$ at $e$, the path $p_e$ contains $s$ (by definition) and leaves $s$ to the South. In particular $p_e$ contains $e^{\pr}$.

Suppose $p_e\neq q$ and let $v$ be the first vertex after $s$ such that $p_e$ and $q$ leave in different directions.  (As before, if $p_e$ and $q$ agree after $s$, then we may apply a similar argument where $p_e$ and $q$ enter some vertex $v$ before $s$ in different directions.)  Let $t$ be the segment $q[(e^{\pr})_{\term},v]$.  If $p_e$ leaves $v$ to the South and $q$ to the East, then $t\in A_q$.  But this implies $s\circ t\in\phi(F)$ so $p_e$ must leave $s\circ t$ to the East, a contradiction.  On the other hand, if $p_e$ leaves $v$ to the East and $q$ to the South, then $s\circ t\in\phi(F)$.  As $\phi(F)=\phi(F)^{\downarrow}$, this implies $t\in\phi(F)$.  In particular, there exist segments $t_1,\ldots,t_l$ such that $t=t_1\circ\cdots\circ t_l$ and $t_i\in A_{p_i}$ for some paths $p_i$ in $F$.  Since $q$ is on top at $e^{\pr}$, it must leave $t_1$ to the East.  Since $q$ and $p_2$ are non-kissing, it leaves $t_2$ to the East as well.  Similarly, it leaves $t_3,\ldots,t_l$ to the East, a contradiction.
\end{poof}

\begin{theorem}\label{thm_quotient}
$\GT(\lambda)$ is a lattice quotient of $\Bic(S)$.
\end{theorem}

\begin{poof}
By Claim \ref{claim_quotient}, the fibers of $\eta$ are the equivalence classes of $\Theta$.  By Claim \ref{claim_covers}, the defining relation of $\GT(\lambda)$ coincides with the covering relations of $\Bic(S)/\Theta$.  Therefore, $\GT(\lambda)$ is a well-defined partial order which is isomorphic to $\Bic(S)/\Theta$.
\end{poof}

As $\Bic(S)$ is a congruence-uniform lattice and congruence-uniformity is preserved by lattice quotients, Theorem \ref{thm_main} follows from Theorem \ref{thm_quotient}.

\section{Cambrian Lattices as Grid-Tamari orders}\label{sec_cambrian}

In this section, we recall the definition of a Cambrian lattice (of type A) as a poset of triangulations of a polygon.  We then identify this lattice with the Grid-Tamari order on a double ribbon shape.

Fix $n\in\Nbb$ and let $Q$ be a directed graph whose underlying graph is a path on $n-1$ vertices.  Label the vertices $v_1,\ldots,v_{n-1}$ in order along this path.  We define a polygon $P$ in $\Rbb^2$ with vertices $w_0,\ldots,w_{n+1}$ such that
\begin{itemize}
\item $w_i$ has $x$-coordinate $i$ for all $i$,
\item $w_0$ and $w_{n+1}$ are above the $x$-axis,
\item $w_1$ and $w_n$ are below the $x$-axis, and
\item for $i=2,\ldots,n-1$, $w_i$ is above the $x$-axis if and only if there is a directed edge $v_{i-1}\ra v_i$ in $Q$.
\end{itemize}

The \emph{Cambrian lattice} $\Camb(Q)$ is the set of triangulations of $P$ whose covering relations are of the form $T\lessdot T^{\pr}$ if $T$ and $T^{\pr}$ differ by a single diagonal and the slope of the diagonal in $T-T^{\pr}$ is less than the slope of the diagonal in $T^{\pr}-T$.

Let $\lambda$ be the double ribbon shape with interior vertices $v_1,\ldots,v_{n-1}$ such that for $i\in\{2,\ldots,n-1\}$,
\begin{itemize}
\item $v_{i-1}$ is North of $v_i$ if $v_{i-1}\ra v_i$ in $Q$ and 
\item $v_{i-1}$ is West of $v_i$ if $v_i\ra v_{i-1}$ in $Q$.
\end{itemize}

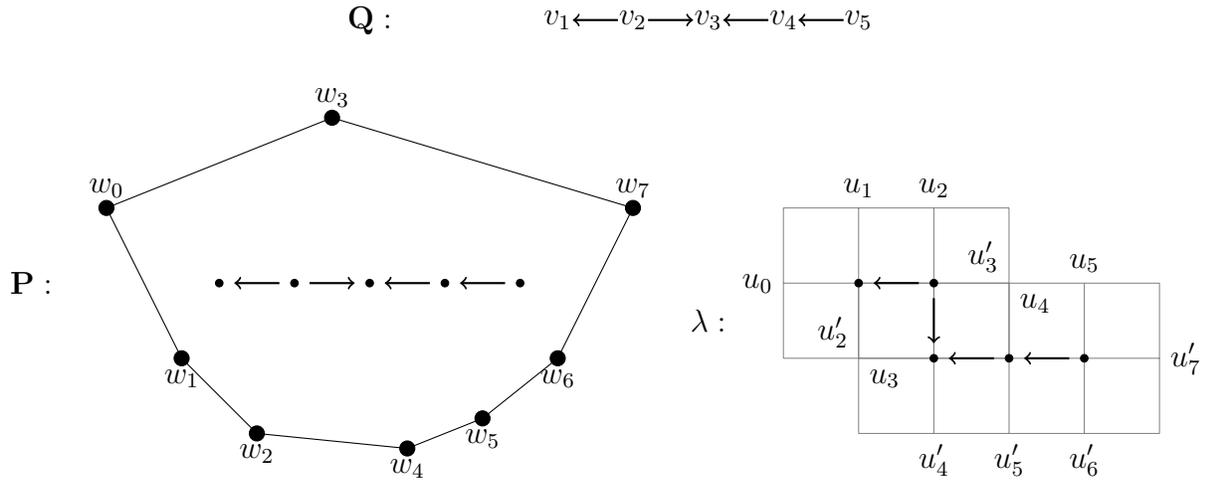
\begin{figure}
\begin{centering}
\begin{tikzpicture}

\begin{scope}[xshift=6cm]
\draw (-2.5,0) node{$\bf Q:$};
\draw[->, thick] (.8,0) -- (.2,0);
\draw[->, thick] (1.2,0) -- (1.8,0);
\draw[->, thick] (2.8,0) -- (2.2,0);
\draw[->, thick] (3.8,0) -- (3.2,0);
\draw (0,0) node{$v_1$};
\draw (1,0) node{$v_2$};
\draw (2,0) node{$v_3$};
\draw (3,0) node{$v_4$};
\draw (4,0) node{$v_5$};
\end{scope}

\begin{scope}[yshift=-3.5cm]
\draw (-1,0) node{$\bf P:$};
\draw (0,1) -- (1,-1) -- (2,-2) -- (4,-2.2) -- (5,-1.8) -- (6,-1) -- (7,1) -- (3,2.2) -- cycle;
\draw (0,1) node[anchor=south]{$w_0$};
\draw (1,-1) node[anchor=north]{$w_1$};
\draw (2,-2) node[anchor=north]{$w_2$};
\draw (3,2.2) node[anchor=south]{$w_3$};
\draw (4,-2.2) node[anchor=north]{$w_4$};
\draw (5,-1.8) node[anchor=north]{$w_5$};
\draw (6,-1) node[anchor=north]{$w_6$};
\draw (7,1) node[anchor=south]{$w_7$};

\filldraw (0,1) circle(1mm);
\filldraw (1,-1) circle(1mm);
\filldraw (2,-2) circle(1mm);
\filldraw (3,2.2) circle(1mm);
\filldraw (4,-2.2) circle(1mm);
\filldraw (5,-1.8) circle(1mm);
\filldraw (6,-1) circle(1mm);
\filldraw (7,1) circle(1mm);

\begin{scope}[xshift=1.5cm]
\draw[->, thick] (.8,0) -- (.2,0);
\draw[->, thick] (1.2,0) -- (1.8,0);
\draw[->, thick] (2.8,0) -- (2.2,0);
\draw[->, thick] (3.8,0) -- (3.2,0);
\filldraw (0,0) circle(.5mm);
\filldraw (1,0) circle(.5mm);
\filldraw (2,0) circle(.5mm);
\filldraw (3,0) circle(.5mm);
\filldraw (4,0) circle(.5mm);
\end{scope}

\end{scope}

\begin{scope}[yshift=-2.5cm, xshift=9cm]

\draw (-1,-1.5) node{$\mathbf \lambda:$};

\draw[step=1cm, gray, very thin] (0,-2) grid (3,0);
\draw[step=1cm, gray, very thin] (1,-3) grid (5,-1);

\draw (0,-1) node[anchor=east]{$u_0$};
\draw (1,0) node[anchor=south]{$u_1$};
\draw (2,0) node[anchor=south]{$u_2$};
\draw (1,-2) node[anchor=south east]{$u_2^{\pr}$};
\draw (1,-2) node[anchor=north west]{$u_3$};
\draw (2,-3) node[anchor=north]{$u_4^{\pr}$};
\draw (3,-1) node[anchor=south east]{$u_3^{\pr}$};
\draw (3,-1) node[anchor=north west]{$u_4$};
\draw (3,-3) node[anchor=north]{$u_5^{\pr}$};
\draw (4,-1) node[anchor=south]{$u_5$};
\draw (4,-3) node[anchor=north]{$u_6^{\pr}$};
\draw (5,-2) node[anchor=west]{$u_7^{\pr}$};

\draw[->, thick] (1.8,-1) -- (1.2,-1);
\draw[->, thick] (2,-1.2) -- (2,-1.8);
\draw[->, thick] (2.8,-2) -- (2.2,-2);
\draw[->, thick] (3.8,-2) -- (3.2,-2);
\filldraw (1,-1) circle(.5mm);
\filldraw (2,-1) circle(.5mm);
\filldraw (2,-2) circle(.5mm);
\filldraw (3,-2) circle(.5mm);
\filldraw (4,-2) circle(.5mm);

\end{scope}

\end{tikzpicture}
\caption{\scriptsize\label{fig_cambrian_quiver}$\GT(\lambda)$ is a Cambrian lattice when $\lambda$ is a double ribbon shape.}
\end{centering}
\end{figure}

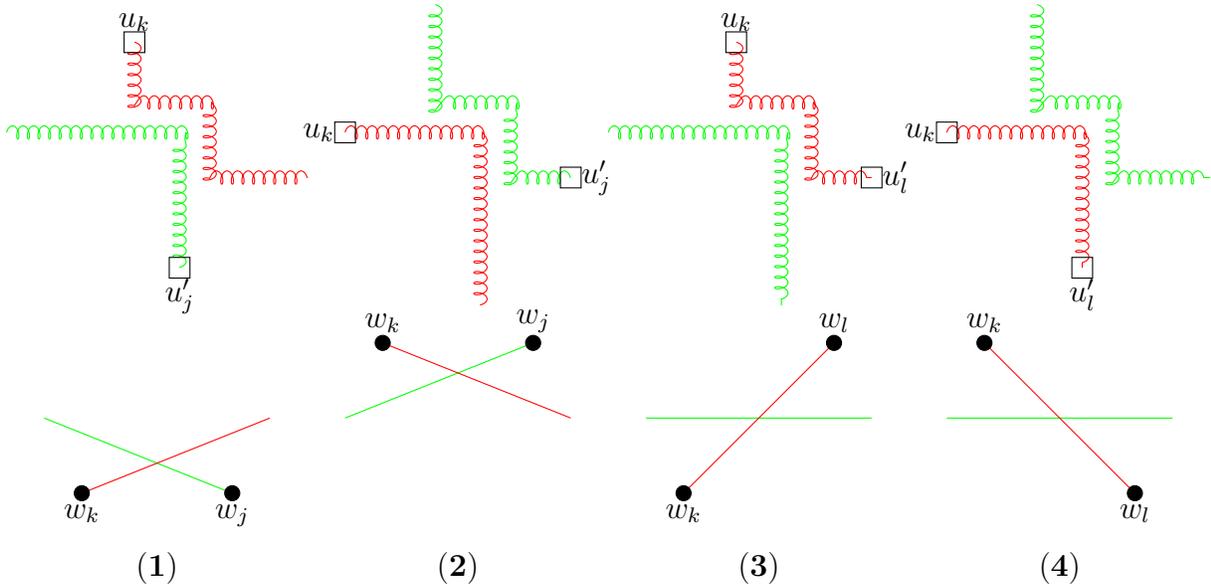
\begin{figure}
\begin{centering}
\begin{tikzpicture}

\begin{scope}

\draw[decorate,decoration={coil,segment length=4pt},color=red] (.2,1) -- (.2,.2) -- (1.2,.2) -- (1.2,-.8) -- (2.5,-.8);
\draw[decorate,decoration={coil,segment length=4pt},color=green] (-1.5,-.2) -- (.8,-.2) -- (.8,-2);

\draw (.2,1) node[minimum size=2mm,draw]{};
\draw (.2,1) node[anchor=south]{$u_k$};
\draw (.8,-2) node[minimum size=2mm,draw]{};
\draw (.8,-2) node[anchor=north]{$u_j^{\pr}$};

\end{scope}

\begin{scope}[xshift=4cm]

\draw[decorate,decoration={coil,segment length=4pt},color=green] (.2,1.5) -- (.2,.2) -- (1.2,.2) -- (1.2,-.8) -- (2,-.8);
\draw[decorate,decoration={coil,segment length=4pt},color=red] (-1,-.2) -- (.8,-.2) -- (.8,-2.5);

\draw (-1,-.2) node[minimum size=2mm,draw]{};
\draw (-1,-.2) node[anchor=east]{$u_k$};
\draw (2,-.8) node[minimum size=2mm,draw]{};
\draw (2,-.8) node[anchor=west]{$u_j^{\pr}$};

\end{scope}

\begin{scope}[xshift=8cm]

\draw[decorate,decoration={coil,segment length=4pt},color=red] (.2,1) -- (.2,.2) -- (1.2,.2) -- (1.2,-.8) -- (2,-.8);
\draw[decorate,decoration={coil,segment length=4pt},color=green] (-1.5,-.2) -- (.8,-.2) -- (.8,-2.5);

\draw (.2,1) node[minimum size=2mm,draw]{};
\draw (.2,1) node[anchor=south]{$u_k$};
\draw (2,-.8) node[minimum size=2mm,draw]{};
\draw (2,-.8) node[anchor=west]{$u_l^{\pr}$};

\end{scope}

\begin{scope}[xshift=12cm]

\draw[decorate,decoration={coil,segment length=4pt},color=green] (.2,1.5) -- (.2,.2) -- (1.2,.2) -- (1.2,-.8) -- (2.5,-.8);
\draw[decorate,decoration={coil,segment length=4pt},color=red] (-1,-.2) -- (.8,-.2) -- (.8,-2);

\draw (-1,-.2) node[minimum size=2mm,draw]{};
\draw (-1,-.2) node[anchor=east]{$u_k$};
\draw (.8,-2) node[minimum size=2mm,draw]{};
\draw (.8,-2) node[anchor=north]{$u_l^{\pr}$};

\end{scope}

\begin{scope}[yshift=-4cm,xshift=.5cm]
\draw[green] (-1.5,0) -- (1,-1);
\draw[red] (-1,-1) -- (1.5,0);

\filldraw (-1,-1) circle(1mm);
\draw (-1,-1) node[anchor=north]{$w_k$};
\filldraw (1,-1) circle(1mm);
\draw (1,-1) node[anchor=north]{$w_j$};

\draw (0,-2) node{$\bf (1)$};

\end{scope}

\begin{scope}[yshift=-4cm,xshift=4.5cm]
\draw[green] (-1.5,0) -- (1,1);
\draw[red] (-1,1) -- (1.5,0);

\filldraw (-1,1) circle(1mm);
\draw (-1,1) node[anchor=south]{$w_k$};
\filldraw (1,1) circle(1mm);
\draw (1,1) node[anchor=south]{$w_j$};

\draw (0,-2) node{$\bf (2)$};

\end{scope}

\begin{scope}[yshift=-4cm,xshift=8.5cm]
\draw[green] (-1.5,0) -- (1.5,0);
\draw[red] (-1,-1) -- (1,1);

\filldraw (-1,-1) circle(1mm);
\draw (-1,-1) node[anchor=north]{$w_k$};
\filldraw (1,1) circle(1mm);
\draw (1,1) node[anchor=south]{$w_l$};

\draw (0,-2) node{$\bf (3)$};

\end{scope}

\begin{scope}[yshift=-4cm,xshift=12.5cm]
\draw[green] (-1.5,0) -- (1.5,0);
\draw[red] (-1,1) -- (1,-1);

\filldraw (-1,1) circle(1mm);
\draw (-1,1) node[anchor=south]{$w_k$};
\filldraw (1,-1) circle(1mm);
\draw (1,-1) node[anchor=north]{$w_l$};

\draw (0,-2) node{$\bf (4)$};

\end{scope}

\end{tikzpicture}
\caption{\label{fig_cambrian_cases}\scriptsize From the proof of Proposition \ref{prop_cambrian}: The four ways two paths may be kissing in $\lambda$ and the corresponding ways two diagonals may cross in $P$.}
\end{centering}
\end{figure}

\begin{proposition}\label{prop_cambrian}
Given $Q$ and $\lambda$ as above, $\Camb(Q)$ is isomorphic to $\GT(\lambda)$ as lattices.
\end{proposition}

\begin{poof}
We first define a bijection between paths in $\lambda$ with diagonals in $P$.  We label the boundary vertices $u_0,\ldots,u_{n-1}$ and $u_2^{\pr},\ldots,u_{n+1}^{\pr}$ where
\begin{itemize}
\item $u_0$ is West of $v_1$ and $u_1$ is North of $v_1$,
\item $u_{n+1}^{\pr}$ is East of $v_{n-1}$ and $u_n^{\pr}$ is South of $v_{n-1}$,
\item for $i\in\{2,\ldots,n-1\}$, if $v_{i-1}$ is North of $v_i$, then $u_i$ is West of $v_i$ and $u_i^{\pr}$ is East of $v_{i-1}$, and
\item for $i\in\{2,\ldots,n-1\}$, if $v_{i-1}$ is West of $v_i$, then $u_i$ is North of $v_i$ and $u_i^{\pr}$ is South of $v_{i-1}$.
\end{itemize}

Every boundary vertex that can start (end) a path is labeled $u_i$ ($u_i^{\pr}$) for a unique $i$.  Let $\tau$ map paths in $\lambda$ with at least one turn to diagonals of $P$ such that the path from $u_i$ to $u_j^{\pr}$ is sent to the diagonal between $w_i$ and $w_j$.  It is straight-forward to check that $\tau$ is a bijection.  We check that two paths $p,p^{\pr}$ are kissing if and only if $\tau(p)$ and $\tau(p^{\pr})$ are crossing.

Let $p$ be the path between $u_i$ and $u_j^{\pr}$, and let $p^{\pr}$ be the path between $u_k$ and $u_l^{\pr}$ for some $i,j,k,l$.  Assume $p$ and $p^{\pr}$ are kissing.  Without loss of generality, we may assume $i<k$.  Then exactly one of the following must hold:
\begin{enumerate}
\item $i<k<j<l$, $u_k$ is North of $v_k$, and $u_j^{\pr}$ is South of $v_{j-1}$;
\item $i<k<j<l$, $u_k$ is West of $v_k$, and $u_j^{\pr}$ is East of $v_{j-1}$;
\item $i<k<l<j$, $u_k$ is North of $v_k$, and $u_l^{\pr}$ is East of $v_{l-1}$; or
\item $i<k<l<j$, $u_k$ is West of $v_k$, and $u_l^{\pr}$ is South of $v_{l-1}$.
\end{enumerate}

Similarly, the diagonal between $w_i$ and $w_j$ crosses the diagonal between $w_k$ and $w_l$ for some $i<j,\ k<l$ in exactly one of the following cases:
\begin{enumerate}
\item $i<k<j<l$ and $w_k$ and $w_j$ are below the $x$-axis;
\item $i<k<j<l$ and $w_k$ and $w_j$ are above the $x$-axis;
\item $i<k<l<j$, $w_k$ is below the $x$-axis, and $w_l$ is above the $x$-axis; or
\item $i<k<l<j$, $w_k$ is above the $x$-axis, and $w_l$ is below the $x$-axis.
\end{enumerate}

Hence, $\tau$ induces an isomorphism of compatibility complexes.  If $F$ and $F^{\pr}$ are adjacent facets of the non-kissing complex, then there exists unique paths $p\in F-F^{\pr}$ and $p^{\pr}\in F-F^{\pr}$.  Checking the four cases above, it is routine to verify that $F<F^{\pr}$ in $\GT(\lambda)$ if and only if the slope of $\tau(p)$ is less than the slope of $\tau(p^{\pr})$.
\end{poof}

\section*{Acknowledgements}

I thank Nathan Reading, Emily Barnard, and Salvatore Stella for their comments about this work.  I also thank Christian Stump for making his data on Grassmann-Tamari orders available for my use.  Lastly, I thank Vic Reiner and my advisor, Pasha Pylyavskyy, for many helpful remarks.

\scriptsize
\bibliography{bib_tamari}{}
\bibliographystyle{plain}

\end{document}